\documentclass[11pt]{amsart}
\usepackage{amssymb}
\begin{document}
\theoremstyle{plain}
\newtheorem{thm}{Theorem}[section]
\newtheorem{theorem}[thm]{Theorem}
\newtheorem{lemma}[thm]{Lemma}
\newtheorem{corollary}[thm]{Corollary}
\newtheorem{proposition}[thm]{Proposition}
\newtheorem{addendum}[thm]{Addendum}
\newtheorem{variant}[thm]{Variant}
\theoremstyle{definition}
\newtheorem{construction}[thm]{Construction}
\newtheorem{notations}[thm]{Notations}
\newtheorem{question}[thm]{Question}
\newtheorem{problem}[thm]{Problem}
\newtheorem{remark}[thm]{Remark}
\newtheorem{remarks}[thm]{Remarks}
\newtheorem{definition}[thm]{Definition}
\newtheorem{claim}[thm]{Claim}
\newtheorem{assumption}[thm]{Assumption}
\newtheorem{assumptions}[thm]{Assumptions}
\newtheorem{properties}[thm]{Properties}
\newtheorem{example}[thm]{Example}
\newtheorem{exercise}[thm]{Exercise}
\newtheorem{conjecture}[thm]{Conjecture}
\newtheorem{criterion}[thm]{Criterion}
\newtheorem{comments}[thm]{Comments}
\newtheorem{blank}[thm]{}
\numberwithin{equation}{thm}
\catcode`\@=11
\def\opn#1#2{\def#1{\mathop{\kern0pt\fam0#2}\nolimits}}
\def\bold#1{{\bf #1}}%
\def\underrightarrow{\mathpalette\underrightarrow@}
\def\underrightarrow@#1#2{\vtop{\ialign{$##$\cr
 \hfil#1#2\hfil\cr\noalign{\nointerlineskip}%
 #1{-}\mkern-6mu\cleaders\hbox{$#1\mkern-2mu{-}\mkern-2mu$}\hfill
 \mkern-6mu{\to}\cr}}}
\let\underarrow\underrightarrow
\def\underleftarrow{\mathpalette\underleftarrow@}
\def\underleftarrow@#1#2{\vtop{\ialign{$##$\cr
 \hfil#1#2\hfil\cr\noalign{\nointerlineskip}#1{\leftarrow}\mkern-6mu
 \cleaders\hbox{$#1\mkern-2mu{-}\mkern-2mu$}\hfill
 \mkern-6mu{-}\cr}}}
\let\amp@rs@nd@\relax
\newdimen\ex@
\ex@.2326ex
\newdimen\bigaw@
\newdimen\minaw@
\minaw@16.08739\ex@
\newdimen\minCDaw@
\minCDaw@2.5pc
\newif\ifCD@
\def\minCDarrowwidth#1{\minCDaw@#1}
\newenvironment{CD}{\@CD}{\@endCD}
\def\@CD{\def\A##1A##2A{\llap{$\vcenter{\hbox
 {$\scriptstyle##1$}}$}\Big\uparrow\rlap{$\vcenter{\hbox{%
$\scriptstyle##2$}}$}&&}%
\def\V##1V##2V{\llap{$\vcenter{\hbox
 {$\scriptstyle##1$}}$}\Big\downarrow\rlap{$\vcenter{\hbox{%
$\scriptstyle##2$}}$}&&}%
\def\={&\hskip.5em\mathrel
 {\vbox{\hrule width\minCDaw@\vskip3\ex@\hrule width
 \minCDaw@}}\hskip.5em&}%
\def\verteq{\Big\Vert&&}%
\def\noarr{&&}%
\def\vspace##1{\noalign{\vskip##1\relax}}\relax\let\amp@rs@nd@&\iffalse}\fi
 \CD@true\vcenter\bgroup\relax\let\\=\cr\iffalse}\fi\tabskip\z@skip\baselineskip20\ex@
 \lineskip3\ex@\lineskiplimit3\ex@\halign\bgroup
 &\hfill$\m@th##$\hfill\cr}
\def\@endCD{\cr\egroup\egroup}
\def\>#1>#2>{\amp@rs@nd@\setbox\z@\hbox{$\scriptstyle
 \;{#1}\;\;$}\setbox\@ne\hbox{$\scriptstyle\;{#2}\;\;$}\setbox\tw@
 \hbox{$#2$}\ifCD@
 \global\bigaw@\minCDaw@\else\global\bigaw@\minaw@\fi
 \ifdim\wd\z@>\bigaw@\global\bigaw@\wd\z@\fi
 \ifdim\wd\@ne>\bigaw@\global\bigaw@\wd\@ne\fi
 \ifCD@\hskip.5em\fi
 \ifdim\wd\tw@>\z@
 \mathrel{\mathop{\hbox to\bigaw@{\rightarrowfill}}\limits^{#1}_{#2}}\else
 \mathrel{\mathop{\hbox to\bigaw@{\rightarrowfill}}\limits^{#1}}\fi
 \ifCD@\hskip.5em\fi\amp@rs@nd@}
\def\<#1<#2<{\amp@rs@nd@\setbox\z@\hbox{$\scriptstyle
 \;\;{#1}\;$}\setbox\@ne\hbox{$\scriptstyle\;\;{#2}\;$}\setbox\tw@
 \hbox{$#2$}\ifCD@
 \global\bigaw@\minCDaw@\else\global\bigaw@\minaw@\fi
 \ifdim\wd\z@>\bigaw@\global\bigaw@\wd\z@\fi
 \ifdim\wd\@ne>\bigaw@\global\bigaw@\wd\@ne\fi
 \ifCD@\hskip.5em\fi
 \ifdim\wd\tw@>\z@
 \mathrel{\mathop{\hbox to\bigaw@{\leftarrowfill}}\limits^{#1}_{#2}}\else
 \mathrel{\mathop{\hbox to\bigaw@{\leftarrowfill}}\limits^{#1}}\fi
 \ifCD@\hskip.5em\fi\amp@rs@nd@}
\newenvironment{CDS}{\@CDS}{\@endCDS}
\def\@CDS{\def\A##1A##2A{\llap{$\vcenter{\hbox
 {$\scriptstyle##1$}}$}\Big\uparrow\rlap{$\vcenter{\hbox{%
$\scriptstyle##2$}}$}&}%
\def\V##1V##2V{\llap{$\vcenter{\hbox
 {$\scriptstyle##1$}}$}\Big\downarrow\rlap{$\vcenter{\hbox{%
$\scriptstyle##2$}}$}&}%
\def\={&\hskip.5em\mathrel
 {\vbox{\hrule width\minCDaw@\vskip3\ex@\hrule width
 \minCDaw@}}\hskip.5em&}
\def\verteq{\Big\Vert&}
\def\novarr{&}
\def\noharr{&&}
\def\SE##1E##2E{\slantedarrow(0,18)(4,-3){##1}{##2}&}
\def\SW##1W##2W{\slantedarrow(24,18)(-4,-3){##1}{##2}&}
\def\NE##1E##2E{\slantedarrow(0,0)(4,3){##1}{##2}&}
\def\NW##1W##2W{\slantedarrow(24,0)(-4,3){##1}{##2}&}
\def\slantedarrow(##1)(##2)##3##4{%
\thinlines\unitlength1pt\lower 6.5pt\hbox{\begin{picture}(24,18)%
\put(##1){\vector(##2){24}}%
\put(0,8){$\scriptstyle##3$}%
\put(20,8){$\scriptstyle##4$}%
\end{picture}}}
\def\vspace##1{\noalign{\vskip##1\relax}}\relax\let\amp@rs@nd@&\iffalse}\fi
 \CD@true\vcenter\bgroup\relax\let\\=\cr\iffalse}\fi\tabskip\z@skip\baselineskip20\ex@
 \lineskip3\ex@\lineskiplimit3\ex@\halign\bgroup
 &\hfill$\m@th##$\hfill\cr}
\def\@endCDS{\cr\egroup\egroup}
\newdimen\TriCDarrw@
\newif\ifTriV@
\newenvironment{TriCDV}{\@TriCDV}{\@endTriCD}
\newenvironment{TriCDA}{\@TriCDA}{\@endTriCD}
\def\@TriCDV{\TriV@true\def\TriCDpos@{6}\@TriCD}
\def\@TriCDA{\TriV@false\def\TriCDpos@{10}\@TriCD}
\def\@TriCD#1#2#3#4#5#6{%
\setbox0\hbox{$\ifTriV@#6\else#1\fi$} \TriCDarrw@=\wd0
\advance\TriCDarrw@ 24pt \advance\TriCDarrw@ -1em
\def\SE##1E##2E{\slantedarrow(0,18)(2,-3){##1}{##2}&}
\def\SW##1W##2W{\slantedarrow(12,18)(-2,-3){##1}{##2}&}
\def\NE##1E##2E{\slantedarrow(0,0)(2,3){##1}{##2}&}
\def\NW##1W##2W{\slantedarrow(12,0)(-2,3){##1}{##2}&}
\def\slantedarrow(##1)(##2)##3##4{\thinlines\unitlength1pt
\lower 6.5pt\hbox{\begin{picture}(12,18)%
\put(##1){\vector(##2){12}}%
\put(-4,\TriCDpos@){$\scriptstyle##3$}%
\put(12,\TriCDpos@){$\scriptstyle##4$}%
\end{picture}}}
\def\={\mathrel {\vbox{\hrule
   width\TriCDarrw@\vskip3\ex@\hrule width
   \TriCDarrw@}}}
\def\>##1>>{\setbox\z@\hbox{$\scriptstyle
 \;{##1}\;\;$}\global\bigaw@\TriCDarrw@
 \ifdim\wd\z@>\bigaw@\global\bigaw@\wd\z@\fi
 \hskip.5em
 \mathrel{\mathop{\hbox to \TriCDarrw@
{\rightarrowfill}}\limits^{##1}}
 \hskip.5em}
\def\<##1<<{\setbox\z@\hbox{$\scriptstyle
 \;{##1}\;\;$}\global\bigaw@\TriCDarrw@
 \ifdim\wd\z@>\bigaw@\global\bigaw@\wd\z@\fi
 \mathrel{\mathop{\hbox to\bigaw@{\leftarrowfill}}\limits^{##1}}
 }
 \CD@true\vcenter\bgroup\relax\let\\=\cr\iffalse}\fi
 \tabskip\z@skip\baselineskip20\ex@
 \lineskip3\ex@\lineskiplimit3\ex@
 \ifTriV@
 \halign\bgroup
 &\hfill$\m@th##$\hfill\cr
#1&\multispan3\hfill$#2$\hfill&#3\\
&#4&#5\\
&&#6\cr\egroup%
\else
 \halign\bgroup
 &\hfill$\m@th##$\hfill\cr
&&#1\\%
&#2&#3\\
#4&\multispan3\hfill$#5$\hfill&#6\cr\egroup \fi}
\def\@endTriCD{\egroup}

\newcounter{Myenumi}
\newenvironment{myenumi}%
{\begin{list}{}{\usecounter{Myenumi}%
\renewcommand{\makelabel}{\arabic{Myenumi}.}%
\settowidth{\leftmargin}{2.n}\settowidth{\labelwidth}{2.n}%
\setlength{\labelsep}{0pt}}}{\end{list}}
\newcounter{Myenumii}
\newenvironment{myenumii}%
{\begin{list}{}{\usecounter{Myenumii}%
\renewcommand{\makelabel}{\alph{Myenumii})}%
\settowidth{\leftmargin}{a)n}\settowidth{\labelwidth}{a)n}%
\setlength{\labelsep}{0pt}}}{\end{list}}
\newcounter{Myenumiii}
\newenvironment{myenumiii}%
{\begin{list}{}{\usecounter{Myenumiii}%
\renewcommand{\makelabel}{\roman{Myenumiii}.}%
\settowidth{\leftmargin}{iv.n}\settowidth{\labelwidth}{iv.n}%
\setlength{\labelsep}{0pt}}}{\end{list}}

\renewenvironment{quote}{\begin{list}{}%
{\settowidth{\leftmargin}{2.n}\setlength{\rightmargin}{0pt}
\renewcommand{\makelabel}{}}%
\item}%
{\end{list}}

\renewenvironment{itemize}%
{\begin{list}{}{\renewcommand{\makelabel}{$\bullet$}%
\settowidth{\leftmargin}{2.n}\settowidth{\labelwidth}{2.n}%
\setlength{\labelsep}{0pt}}}{\end{list}}

\newsymbol\onto 1310
\def\into{\DOTSB\lhook\joinrel\rightarrow}
\newcommand{\sA}{{\mathcal A}}
\newcommand{\sB}{{\mathcal B}}
\newcommand{\sC}{{\mathcal C}}
\newcommand{\sD}{{\mathcal D}}
\newcommand{\sE}{{\mathcal E}}
\newcommand{\sF}{{\mathcal F}}
\newcommand{\sG}{{\mathcal G}}
\newcommand{\sH}{{\mathcal H}}
\newcommand{\sI}{{\mathcal I}}
\newcommand{\sJ}{{\mathcal J}}
\newcommand{\sK}{{\mathcal K}}
\newcommand{\sL}{{\mathcal L}}
\newcommand{\sM}{{\mathcal M}}
\newcommand{\sN}{{\mathcal N}}
\newcommand{\sO}{{\mathcal O}}
\newcommand{\sP}{{\mathcal P}}
\newcommand{\sQ}{{\mathcal Q}}
\newcommand{\sR}{{\mathcal R}}
\newcommand{\sS}{{\mathcal S}}
\newcommand{\sT}{{\mathcal T}}
\newcommand{\sU}{{\mathcal U}}
\newcommand{\sV}{{\mathcal V}}
\newcommand{\sW}{{\mathcal W}}
\newcommand{\sX}{{\mathcal X}}
\newcommand{\sY}{{\mathcal Y}}
\newcommand{\sZ}{{\mathcal Z}}
\newcommand{\ssA}{{\mathfrak A}}
\newcommand{\ssB}{{\mathfrak B}}
\newcommand{\ssC}{{\mathfrak C}}
\newcommand{\ssD}{{\mathfrak D}}
\newcommand{\ssE}{{\mathfrak E}}
\newcommand{\ssF}{{\mathfrak F}}
\newcommand{\ssG}{{\mathfrak G}}
\newcommand{\ssH}{{\mathfrak H}}
\newcommand{\ssI}{{\mathfrak I}}
\newcommand{\ssJ}{{\mathfrak J}}
\newcommand{\ssK}{{\mathfrak K}}
\newcommand{\ssL}{{\mathfrak L}}
\newcommand{\ssM}{{\mathfrak M}}
\newcommand{\ssN}{{\mathfrak N}}
\newcommand{\ssO}{{\mathfrak O}}
\newcommand{\ssP}{{\mathfrak P}}
\newcommand{\ssQ}{{\mathfrak Q}}
\newcommand{\ssR}{{\mathfrak R}}
\newcommand{\ssS}{{\mathfrak S}}
\newcommand{\ssT}{{\mathfrak T}}
\newcommand{\ssU}{{\mathfrak U}}
\newcommand{\ssV}{{\mathfrak V}}
\newcommand{\ssW}{{\mathfrak W}}
\newcommand{\ssX}{{\mathfrak X}}
\newcommand{\ssY}{{\mathfrak Y}}
\newcommand{\ssZ}{{\mathfrak Z}}

\newcommand{\A}{{\mathbb A}}
\newcommand{\B}{{\mathbb B}}
\newcommand{\C}{{\mathbb C}}
\newcommand{\D}{{\mathbb D}}
\newcommand{\E}{{\mathbb E}}
\newcommand{\F}{{\mathbb F}}
\newcommand{\G}{{\mathbb G}}
\newcommand{\HH}{{\mathbb H}}
\newcommand{\I}{{\mathbb I}}
\newcommand{\J}{{\mathbb J}}
\newcommand{\K}{{\mathbb K}}
\renewcommand{\L}{{\mathbb L}}
\newcommand{\M}{{\mathbb M}}
\newcommand{\N}{{\mathbb N}}
\renewcommand{\P}{{\mathbb P}}
\newcommand{\Q}{{\mathbb Q}}
\newcommand{\R}{{\mathbb R}}
\newcommand{\bS}{{\mathbb S}}

\newcommand{\T}{{\mathbb T}}
\newcommand{\U}{{\mathbb U}}
\newcommand{\V}{{\mathbb V}}
\newcommand{\W}{{\mathbb W}}
\newcommand{\X}{{\mathbb X}}
\newcommand{\Y}{{\mathbb Y}}
\newcommand{\Z}{{\mathbb Z}}
\newcommand{\id}{{\rm id}}
\newcommand{\rank}{{\rm rank}}
\newcommand{\END}{{\mathbb E}{\rm nd}}
\newcommand{\End}{{\rm End}}
\newcommand{\Hg}{{\rm Hg}}
\newcommand{\tr}{{\rm tr}}
\newcommand{\Tr}{{\rm Tr}}
\newcommand{\Sl}{{\rm Sl}}
\newcommand{\Gl}{{\rm Gl}}
\newcommand{\Cor}{{\rm Cor}}
\newcommand{\GL}{{\rm GL}}
\newcommand{\MT}{{\rm MT}}
\newcommand{\Hdg}{{\rm Hdg}}
\newcommand{\MTV}{{\rm MTV}}
\newcommand{\SO}{{\rm SO}}
\newcommand{\Hom}{{\rm Hom}}
\newcommand{\Ker}{{\rm Ker}}
\newcommand{\Lie}{{\rm Lie}}
\newcommand{\Aut}{{\rm Aut}}
\newcommand{\Image}{{\rm Image}}
\newcommand{\Gr}{{\rm Gr}}
\newcommand{\Id}{{\rm Id}}
\newcommand{\rk}{{\rm rk}}
\newcommand{\pardeg}{{\rm par.deg}}
\newcommand{\SU}{{\rm SU}}
\newcommand{\Res}{{\rm Res}}
\newcommand{\Sym}{{\rm Sym}}
\title[Rigidity for Families of Polarized Calabi-Yau Varieties]{Rigidity for Families of Polarized Calabi-Yau Varieties}
\author[Yi Zhang ]{Yi Zhang}
\address{Center of Mathematical Sciences\\
Zhejiang University (Mailbox 1511) \\
Hangzhou 310027\\ P.R.China \\} \email{yzhang@cms.zju.edu.cn}
\thanks{The research was partially supported by The Institute of Mathematical Sciences and
the Department of Mathematics at The Chinese University of Hong
Kong.}

\maketitle

\begin{abstract} In this paper, we study the analogue of the Shafarevich conjecture
for polarized Calabi-Yau varieties. We use variations of Hodge
structures and Higgs bundles to establish a criterion for the {\it
rigidity} of families. We then apply the criterion to obtain that
some important and typical families of Calabi-Yau varieties are
rigid, for examples., Lefschetz pencils of Calabi-Yau varieties,
{\it strongly degenerated} families (not only for families of
Calabi-Yau varieties), families of Calabi-Yau varieties admitting
a degeneration with {\it maximal unipotent monodromy}.
\end{abstract}

\setcounter{tocdepth}{1} \setcounter{page}{1}
\setcounter{section}{-1}


\section{Introduction}
Throughout this paper,  the base field $k$ is the complex
number field $\C.$ \\

\subsection*{Shafarevich conjecture over function field}

 At the 1962 ICM in Stockholm,  Shafarevich conjectured:
 {\it \textquotedblleft There exists only a finite
number of fields of algebraic functions }$K/k$ {\it of a given
genus }$g \geq1,$ {\it the critical prime divisors of which belong
to a given finite set }$S$ " (cf. \cite{Sh}). Shafarevich proved
his own conjecture in the setting of hyperelliptic curves in one
unpublished work.  The conjecture was confirmed by Par\v sin for
$S=\emptyset,$ by Arakelov in general.

Let $C$ be a smooth projective curve of genus $g(C)$ and $S
\subset C$ be a finite subset. An algebraic family over $C$ is
called {\it isotrivial} if over an open dense subscheme of $C$ any
two smooth fibers are isomorphic. The geometric description of the
Shafarevich conjecture is:
\begin{center} {\it Fix $(C, S)$ and an integer $q\geq 2,$ there
are only finitely many non-isotrivial smooth families of curves of
genus $q$ over $C\setminus S.$}
\end{center}

Let $\sM_q$ be the coarse moduli space of smooth projective curves
of genus $q\geq 2.$ Adding stable curves at the boundary of
$\sM_q,$ one has the Deligne-Mumford compactification
$\overline{\sM_q}.$ Fix $(C,S)$ and $q\geq 2,$ a smooth family
$(f:X\rightarrow C\setminus S)$ of curves of genus $q$ induces
naturally a unique moduli morphism $\eta_f:C\setminus S
\rightarrow \sM_q.$  $\eta_f$ can be extended to
$\overline{\eta}_f: C\rightarrow \overline{\sM}_q$ because of the
smoothness of $C.$ Hence, parameterizing families is same as
parameterizing these morphisms which can be characterized by their
graphs. The graph $\Gamma_{\overline{\eta}_f}$ is a projective
curve contained in $C\times \overline{\sM}_q$ with the first
projection mapping itself isomorphically onto $C ,$ hence the
problem is translated into looking for a parametrization $\T$ in
the Hilbert scheme of $C\times \overline{\sM}_q.$  In general, the
Hilbert scheme is an infinite union of schemes of finite type,
each component has a different Hilbert polynomial and represents a
deformation type of family. Fortunately, this parameterizing
scheme $\T$ is of finite type over $\C,$ i.e., only finite Hilbert
polynomials can actually occur. The original conjecture was
reformulated by Arakelov and Par\v sin into four problems (cf.
\cite{Ara}\cite{Par}):
\begin{conjecture}[Shafarevich Conjecture] Fix $(C,S)$ and $q\geq 2.$
\begin{enumerate}
    \item [({\bf B})] The elements of the set of non-isotrivial families of
curves of genus $q$ over $C$ with singular locus $S$ are
parameterized by points of a scheme $\T $ of finite type over $\C$
({\it Boundedness}).
    \item [({\bf R})] Any deformation of a
non-isotrivial family of curves of genus $q$ over $C$ with
singular locus $S$ is trivial, i.e., $\dim \T=0$ ({\it Rigidity}).
    \item [({\bf H})] No non-isotrivial family of
curves of genus $q$ exists if $2g(C)-2+\#S \leq 0,$ i.e., $\T \neq
\emptyset \Rightarrow 2g(C)-2+\#S > 0 $ ({\it Hyperbolicity}).
    \item [({\bf WB})] For a non-isotrivial
family $f:X\rightarrow C ,$ $\deg f_*\omega_{X/C}^m$ is bounded
above in term only of $ g(C),\#S,q,m$ ({\it Weak Boundedness}).
\end{enumerate}
\end{conjecture}
The Hilbert polynomial of $\Gamma_{\overline{\eta}_f}$ is
determined by its first term $\deg \overline{\eta}_f^*\sL$ where
$\sL$ is a fixed ample line on $\overline{\sM_q}.$  Thus, ({\bf
B}) is equivalent to the boundedness of
$\deg\overline{\eta}_f^*\sL$ due to Mumford's works on moduli
spaces of curves, i.e., that it is sufficient to prove ({\bf WB}).
On the other hand, ({\bf WB}) is obtained by the well-known
Arakelov inequality.  It was shown recently that $({\bf WB})
\Rightarrow ({\bf H})$ if the general fiber is a smooth curve of
$q\geq 2$ (cf. \cite{Kovacs}). ({\bf R}) follows directly from the
positivity of relative dual sheafs of non-isotrivial
families.\\

\subsection*{Shafarevich problems for higher dimensional varieties}
Define $Sh(C, S, K)$ to be the set of all equivalent classes of
non-isotrivial family $\{ f:\sX \rightarrow C \}$ such that
$\sX_b$ is a smooth projective variety with type $'K'$ for any
$b\in C\setminus S.$  Two such families are equivalent if they are
isomorphic over $C\setminus S.$ The general Shafarevich problem is
to find $'K'$ and the data $(C, S)$ such that the set $Sh (C, S,
K)$ is finite.

\begin{example}
Faltings dealt with the case of Abelian varieties, and he
formulated a Hodge theoretic condition (the Deligne-Faltings
$(*)$-condition)  for a fiber space to be rigid (cf. \cite{Falt}):

{\it A smooth family $f:\sX_0 \rightarrow C\setminus S$ of Abelian
varieties is said to satisfy $(*)$-condition if that any
anti-symmetric endomorphism $\sigma$ of $\V_\Z=R^1f_*\Z$ defines an
endomorphism of $\sX_0$ (so $\sigma$ is of type $(0, 0)$).}

By the global Torelli theorem, a polarization of the Abelian scheme
induces a sympletic bilinear form $Q$ on $R^1f_*\Z.$ Faltings showed
that the $(*)$-condition is equivalent to
$$\End^Q(\V_\Q)\otimes \C=(\End^Q(\V_\Q)\otimes \C)^{0,0}. $$
On the other hand, the Zariski tangent space of the moduli space
of Abelian schemes over $C-S$ with a fixed polarization is
isomorphic to $$(\End^Q(\V_\Q)\otimes \C)^{-1,1}.$$
\end{example}

\begin{example}If a family is not rigid, one should have nonrigid VHSs.
By utilizing the differential geometry of period maps and Hodge
metrics on period domains, Peters generalized the result of
Faltings to polarized variations of Hodge structure of arbitrary
weight (cf. \cite{Pe}): {\it A polarized variation of Hodge
structure $(\V,Q)$ underlying $\V_\Z$ is rigid if and only if
$$(\End^Q(\V_\Q)\otimes \C)^{-1, 1}=0. $$ }
\end{example}
\begin{example}As studying deformations of a family can be reduced to
studying deformations of the corresponding period map, Jost and
Yau analyzed $Sh(C,S,K)$ for a large class of varieties by
harmonic maps (cf. \cite{JY}). They provided analytic methods to
solve the {\it rigidity} and gave differential geometric proofs of
Shafarevich conjectures. They started to study Higgs bundles with
singular Hermitian metrics and their applications to Shafarevich
problems.
\end{example}
Eyssidieux and Mok also have many results on deformations of
period maps in case that period domains are Hermitian symmetric
(cf. \cite{Mok3},\cite{E1,E2},\cite{EM}). They showed the {\it gap
rigidity} (which is stronger than the {\it rigidity} in our case)
for locally bounded symmetric domains of certain type (including
all tube domains).
\begin{example}[Nonrigid Family] Faltings constructed an example to show that
the set $Sh(C,S,K)$ is infinite for Abelian varieties of dimension
$\geq 8$ with some type $'K'$ (cf. \cite{Falt}). Saito and Zucker
generalized the construction of Faltings to the setting in case
that the condition $'K'$ is an algebraic polarized $K3$ surface,
and they were able to classify all cases if the set $Sh (C,S,K)$
is infinite (cf. \cite{SZ}).
\end{example}

All these examples did not require that families are polarized.
For curves, the condition $q\geq 2$ is equivalent to that the
canonical line bundle $\omega_{\sX_{gen}}$ of a general fiber is
ample, hence that fixing $g(\sX_{gen})$ can be replaced by that
fixing the Hilbert polynomial $h_{\omega_{\sX_{gen}}}$ of
$\omega_{\sX_{gen}},$ and these families automatically become
families of canonical polarized curves. 
For higher dimensional varieties, instead of fixing $g(\sX_{gen})$
we should consider any polarized projective variety $(X,L)$ such
that $L$ is an ample line bundle on $X$ with $\chi(X,L^{\nu})
\equiv \mbox{ a fixed Hilbert polynomial }h(\nu)$ for $\nu >>0.$

Fix a pair $(C,S)$ such that $C$ is a nonsingular projective curve
and $S$ is a set of finite points of $C.$ Consider the
polarization, $Sh(C,S,K)$  is now defined to be the set of all
equivalent classes of non-isotrivial smooth polarized family $\{(
f:\sX \to C\setminus S, \sL)\}$ satisfying: $\sX_b=f^{-1}(b)$ is a
smooth projective variety with type $'K';$ $\sL$ is an invertible
sheaf on $\sX,$ relatively ample over $C\setminus S$ with a fixed
Hilbert polynomial $\chi(\sL^\nu|_{\sX_b})\equiv h(\nu)$ for any
$b\in C\setminus S$ and $\nu >>0;$ two families are equivalent if
they are isomorphic over $C\setminus S$ as polarized families (see
the following theorem \ref{MGV}). Therefore, the analogue of the
Shafarevich conjecture for higher dimensional varieties is
formulated as follows:
\begin{conjecture}[The Shafarevich Problem for Higher
Dimensional Polarized Varieties] Fix a pair $(C,S)$ and a Hilbert
polynomial $h(\nu).$
\begin{enumerate}
    \item [({\bf B})] The elements of $Sh(C,S,K)$ are parameterized by
points of a scheme $\T$ of finite type over $\C.$
    \item [(\bf{R})]  When does one have $\dim \T=0?$
    \item [(\bf{H})]  $\T\neq \emptyset \Rightarrow 2g(C)-2+\#S>0.$
    \item [(\bf{WB})] For
a family $(f:X\rightarrow C) \in Sh(C,S,K),$  $\deg
f_*\omega_{X/C}^m$ is bounded above in term of $g(C), \#S, h, m. $
In particular,  the bound is independent of $f. $
\end{enumerate}
\end{conjecture}

{\bf Remark.} ${\bf (WB)\Rightarrow (H)}$ holds for any
canonically polarized
family. 
\begin{example}(Recent Results).
\begin{myenumi}
    \item Migliorni, Kov\'ac and Zhang proved that families of minimal
algebraic surfaces of general type over a curve of genus $g$ and
$\#S$ singular fibers such that $2g(C)-2+\#S\leq0 $ are isotrivial
(cf. \cite{Mi}\cite{Ko}\cite{Z}). Oguiso and Viehweg proved it for
families of elliptic surfaces (cf. \cite{O-V}).
    \item Jost-Zuo and Viehweg-Zuo recently made contributions
    to the {\it weak boundedness} (cf. \cite{JZ3}\cite{VZ6}), and obtained Arakelov-Yau type inequalities for
    families of higher dimension varieties over a curve.
    \item Consider an algebraic family over a fixed smooth curve, let $F$ be a general fiber.
Bedulev-Viehweg proved {\bf (B)} if $F$ is an algebraic surface of
general type. They also proved {\bf (WB)} if $F$ is a canonically
polarized variety (cf. \cite{BV}). As explained in \cite{VZ2} one
actually obtains {\bf (B)}. Precisely, Viehweg-Zuo obtained that
${\bf (WB)} \Rightarrow {\bf (B)}$ if $F$ is a minimal model of
general type or if the family is not degenerate and $F$ is a
minimal model of Kodaira dimension zero (cf. \cite{VZ2,VZ4}), and
they also showed that {\bf (WB)} holds if $F$ is a minimal model
of general type or $F$ has  semi-ample $\omega_F.$
    \item Liu-Todorov-Yau-Zuo gave another proof of the {\it boundedness} of the analogue of
Shafarevich conjecture for Calabi-Yau manifolds by originally
using the {\it Schwarz-Yau lemma} and the {\it Bishop
compactness}. They also constructed a nonrigid family of
Calabi-Yau manifolds and related {\it Yukawa couplings} with the
{\it rigidity} (\cite{LTYZ}).
\end{myenumi}
\end{example}
Moduli spaces of Calabi-Yau manifolds play a pivotal role in the
classification theory of Calabi-Yau varieties and in Mirror
Symmetry. Unfortunately, we know little about their structures. As
studying moduli stacks can be reduced to studying families of
manifolds, we begin to study the analogue of the Shafarevich
conjecture for families of high dimensional polarized varieties
and understand part of the structure of moduli stacks. The
existence of coarse moduli spaces of polarized manifolds was
proven by Mumford, Gieseker and Viehweg, it is the fundamental
theorem for us to study moduli problems for manifolds.
\begin{theorem}[cf. \cite{V95}]\label{MGV}
Let $h$ be a fixed polynomial of degree $n$ with $h(\Z)\subset
\Z.$ Define moduli functor

$ \sM_h(Y):= \{ (f: \sX \rightarrow Y, \sL); f \mbox{ flat,projective and $\sL$ invertible,relatively}\\
\mbox{ample over $Y,$ such that: for all $p\in Y(\C)$ $\sX_p=f^{-1}(p)$ is a projective}\\
\mbox{manifold with semi-ample canonical bundle and }
\chi(\sL|_{\sX_p})=h\, \}/\sim.$

Then, the moduli functor $\sM_h$ is bounded by the Matsusaka Big
theorem, and there exists a quasi-projective coarse moduli scheme
$M_h$ for $\sM_h,$ of finite type over $\C.$  Moreover, if
$\omega_\Gamma^\delta=\sO_{\Gamma} \forall \,\Gamma\in \sM_h(\C)$
for one integer $\delta > 0,$  then for some $p>0$ there exists an
ample line bundle $\lambda^{(p)}$ on $M_h$ such that
$\phi_g^*\lambda^{(p)}=g_*\omega_{\sX/Y}^{\delta\cdot p} $ for any
family $(g: \sX \to Y, \sL)\in \sM_h(Y)$ with moduli morphism
$\phi_g: Y \to M_h.$
\end{theorem}

{\bf Remarks.} A line bundle $\sB$ is semi-ample if for some
$\mu>0$ the sheaf $\sB^{\mu}$ is generated by global sections.
$(f: \sX \to Y, \sL)\sim (f': \sX'\to Y, \sL')$ if there exist a
$Y$-isomorphism $\tau: \sX \to \sX'$ and an invertible sheaf $\sF$
on $Y$ such that $\tau^*\sL'\cong \sL\otimes f^*\sF.$ The
statement is also true if \textquotedblleft $\cong$" is replaced
by the numerical equivalence \textquotedblleft $\equiv$".\\

Recently, Viehweg-Zuo obtained remarkable results for the
Shafarevich problem.
\begin{theorem}[cf. \cite{VZ1,VZ2,VZ3}]
Brody hyperbolicity holds for moduli spaces of canonically
polarized complex manifolds, thus the {\it boundedness} of
$Sh(C,S,Z)$ for arbitrary $Z$ with $\omega_{Z}$ semi-ample holds.
Moreover, the automorphism group of moduli stacks of polarized
manifolds is finite and the rigidity holds for a general family.
\end{theorem}
{\bf Remark.} A complex analytic space $\sN$ is called {\it Brody
hyperbolic} if every holomorphic map $\C \to \sN$ is constant. Its
algebraic version is called {\it algebraic hyperbolic}. One has
that ({\bf H}) is true if the moduli space is {\it algebraic
hyperbolic}.

\subsection*{Rigidity for families of Calabi-Yau varieties}

In this paper, a {\it Calabi-Yau manifold} is a smooth projective
variety $X$ (of dimension $n$) such that the canonical line bundle
$K_X$ is trivial and $h^{i}(X,\sO_{X})=0 \,\forall i$ with
$0<i<n.$ We study the {\it rigidity} problem of the analogue
Shafarevich conjecture for Calabi-Yau manifolds. Before we cite
our main results, we shall point out that the {\it rigidity} for
the analogue of the Shafarevich conjecture fails for general
condition by the following key observation.
\begin{example}(cf. \cite{VZ7})
Very recently, Viehweg-Zuo constructed a nontrivial family of
Calabi-Yau manifolds such that the closed fibers are Calabi-Yau
manifolds endowed with {\it complex multiplication} over a dense
set. Precisely, they obtained a family $g:\sZ \to S$ of quintic
hypersurfaces in $\C\P^4$ such that $S$ is finite dominant over a
ball quotient (which is a Shimura variety) and $S$ has a dense set
of {\it CM} points. Furthermore, they got an important
counterexample for the {\it rigidity} part of the Shafarevich
problem by showing that there exists a product of moduli spaces of
hypersurfaces of degree $d$ in $\P^n$ and that this product can be
embedded into the moduli space of hypersurfaces of degree $d$ in
$\P^N$ for one $N>n.$
\end{example}
Altogether, the remaining step for the Shafarevich problem is to
find nice conditions for rigidity. We obtain three main results as
follows:
\begin{enumerate}
    \item[(I)] We prove that any non-isotrivial Lefschetz pencil
$f: \ssX\to \P^1$ of Calabi-Yau varieties of odd dimension $n$ is
rigid (Lefschetz pencils of even dimensional Calabi-Yau varieties
are automatically trivial). The proof depends on the construction
of Lefschetz pencils in Deligne's Weil conjecture I. Let $f^0$ be
the maximal smooth subfamily of $f$ and $\V$ be the $\Z$-local
system of vanishing cycles space. Actually, we obtain that the
pieces of $(n,0)$-type and $(0,n)$-type of the VHS $R^nf^0_*(\C)$
are both in $\V_\C=\V\otimes \C.$ On the other hand, if the family
$f$ is nonrigid we would have a nonzero $(-1,1)$-type endomorphism
$\sigma$ of $\V_\C$ which is flat under the Gauss-Manin
connection. The action of $\sigma$ induces a nontrivial splitting
of the local system $\V_\C,$ but it contradicts to that $\V$ is
absolutely irreducible.


    \item[(II)]  We obtain a general result that any non-isotrivial {\it strongly
degenerated} family (not only for families of Calabi-Yau
varieties) must be rigid. A family over projective smooth curve is
called {\it strongly degenerate} if it has a singular fiber with
only pure Hodge type cohomology. If the family is nonrigid we then
have a nonzero $(-1,1)$-type endomorphism $\sigma $ which is flat
under the Gauss-Manin connection as same as in (I). The K\"unneth
formula says that we can identify this $\sigma$ with a
monodromy-invariant section of a VHS from the self-product family.
We compare the Hodge type of $\sigma$ with the Hodge types of
cohomology groups of the singular fiber of the self-product
family. Then, we have a contradiction that $\sigma$ must be zero.
As a corollary, we obtain the {\it Weakly Arakelov theorem} for
high dimensional varieties.


    \item [(III)] We introduce a general criterion of Viehweg-Zuo for
    rigidity. From the criterion, we deduce a result of Liu-Todorov-Yau-Zuo and
Viehweg-Zuo: a family of Calabi-Yau varieties over an algebraic
curve is rigid if its {\it Yukawa coupling} is nonzero. Together
with the results of Schmid and Simpson on residues of holomorphic
vector bundles over singularizes, we prove that any family of
Calabi-Yau varieties over an algebraic curve admitting a
degeneration with {\it maximal unipotent monodromy} must be rigid.
\end{enumerate}

\section{Higgs Bundles over Quasi-Projective Manifolds}

\subsection*{The generalized Donaldson-Simpson-Uhlenbeck-Yau correspondence}

Let the base $M$ be a quasi-projective manifold such that there is
a smooth projective completion $\overline{M}$ with a reduced
normal crossing divisor $D_{\infty}=\overline{M}-M.$\\

Let $(V,\nabla)$ be a flat $\GL(n,\C)$ vector bundle on $M,$ i.e.,
a fundamental representation $\rho :\pi_1(M)\rightarrow
\GL(n,\C).$ A Hermitian metric $H$ on $V$ leads to a decomposition
$\nabla=D_H+\vartheta,$ which corresponds to the Cartan
decomposition of Lie algebra $gl(n,\C)= \mathfrak u(n)\oplus
\mathfrak p. $ $D_H$ is a unitary connection preserving the metric
$H,$ so $\rho$ is unitary if and only if $\vartheta =0.$ The
Hermitian metric $H$ can be regarded as a $\rho$-equivariant map
$$h:  \widetilde{M} \rightarrow \GL(n,\C)/\mathrm{U}(n) \mbox{ with }dh=\vartheta$$
where $\widetilde{M}$ is the universal covering of $M.$ With
respect to the complex structure of $M,$ one has the decomposition
$$D_H=D_H^{1,0}+D_H^{0,1}, \,
            \vartheta=\vartheta^{1,0}+\vartheta^{0,1}.$$
The following three conditions are equivalent:
\begin{itemize}
    \item $\nabla^*_H(\vartheta)=0 \mbox{ ($\nabla^*_H$ is defined by $(e,\nabla^*_H(f))_H:= (\nabla(e),f)_H$)};$
    \item $(D_H^{0,1})^2=0, \  D_H^{0,1}(\vartheta^{1,0})=0, \
         \vartheta^{1,0}\wedge \vartheta^{1,0}=0;$
    \item $h$ is a harmonic map, i.e., $H$ is a {\it harmonic metric} (or
          $V$ is {\it harmonic}).
\end{itemize}
Altogether, suppose that $H$ is {\it harmonic},
$(E,\overline{\partial}_E,\theta)$ is a Higgs bundle with respect
to the holomorphic structure $\overline{\partial}_E:=D_H^{0,1}$
where $E$ takes the underlying bundle as $V$ and
$\theta:=\vartheta^{1,0};$ $D_H$ is the unique metric connection
with respect to $\overline{\partial}_E;$ and $H$ is the {\it
Hermitian-Yang-Mills metric} on the Higgs bundle
$(E,\overline{\partial}_E,\theta),$ i.e., $$D_H^2=-(
\theta_H^*\wedge\theta + \theta\wedge\theta_H^*).$$ The existence
of the $\rho$-equivariant harmonic map was proven by Simpson in
case that $M$ is a curve (cf. \cite{Si2}), by Jost-Zuo in case
that $M$ is a higher dimensional manifold (cf. \cite{JZ1,JZ2}). In
fact, Jost-Zuo obtained some useful estimates of curvatures at the
infinity as  Simpson did for punctured curves (cf. Section 2, Main
Estimate, Theorem 1 in \cite{Si2}). The $\rho$-equivariant
harmonic map is unique and depends only on $\rho$ if $M$ is
compact, but the uniqueness does not hold if $M$ is not compact.
On the other hand, a theorem of Cornalba and Griffiths for
extensions of analytic sheaves (cf. \cite{CG}) shows that the
induced Higgs bundle $E$ can extend to a coherent sheaf
$\overline{E}$ on $\overline{M}$ and $\theta$ can extend to
$\overline{\theta}\in \Gamma(\overline{M}, \sE
nd(\overline{E})\otimes \Omega^1_{\overline{M}}(\log
D_{\infty})).$ The extension of $(E,\theta)$ is not unique, but
one can treat this nonuniqueness by taking filtered extensions
$(E,\theta)_{\alpha},$ and obtains a filtered Higgs bundle
$\{(E,\theta)_{\alpha} \}$
(cf. \cite{Zuo1}).\\

Conversely, let $(E,\overline{\partial}_E,\theta)$ be a Higgs
bundle equipped with a Hermitian metric  $H.$ One has a unique
metric connection $D_H$ on $(E,\theta)$ with respect to the
holomorphic structure $\overline{\partial}_E,$ and a $(0,1)$-form
$\theta^*_H$ is determined by $(\theta e,f)_H=(e,\theta^*_H f)_H$
where $e , f$ are arbitrary sections of $E.$ Denote
$\partial_E:=D_H-\overline{\partial}_E,$  $\nabla^{'}:=
\partial_E +\theta^*_H,$ and  $\nabla^{''}:=\overline{\partial}_E+\theta.$
One has $\nabla^{''}\circ \nabla^{''}=0,$ and $(\partial_E)^2=0 $
as $D_H^2=\pi^{1,1}(D_H^2).$ Denote $\nabla=
\nabla^{'}+\nabla^{''}.$ Then, $\nabla=D_H+ \theta+\theta^*_H.$
Let $(V,\nabla^{''} )$ be another holomorphic vector bundle with
respect to $\nabla^{''}$ where $V$ takes the underlying bundle as
$E.$ The metric $H$ on $(E,\theta)$ is called {\it
Hermitian-Yang-Mills} if $\nabla$ is a flat connection on $V,$
i.e., $H$ is a {\it harmonic metric} on $V$ (cf.
\cite{JY},\cite{Zuo1}). If $M$ is a compact K\"ahler manifold or a
quasi-projective curve, one has the {\it
Donaldson-Simpson-Uhlenbeck-Yau correspondence} (DSUY
correspondence), i.e., that the {\it Hermitian-Yang-Mills metric}
exists (cf. \cite{Do1},\cite{Si1},\cite{Si2},\cite{UY}). Suppose
that $(V,\nabla)$ is flat, as
$(\nabla^{1,0})^2=\pi^{2,0}(\nabla^2)=0 $ (with respect to
$\overline{\partial}_E$), then
$\partial_E(\theta)=\overline{\partial}_E(\theta_H^*)=0;$ hence
$\nabla^{'}\circ \nabla^{'}=0$ and $\nabla^{'}$ is just the
Gauss-Manin connection of the holomorphic bundle
$(V,\nabla^{''}).$
\begin{corollary}
Let $(E,\overline{\partial}_E,\theta,H)$ be a Higgs bundle induced
from a {\it harmonic bundle}. Then, one has 
\begin{equation}
\begin{tabular}{|c|c|c|}
 \hline
  $\overline{\partial}_E(\theta)=0$ & $\partial_E(\theta)=0$ & $\theta\wedge\theta=0$ \\
\hline
  $\overline{\partial}_E(\theta_H^*)=0$ & $\partial_E(\theta_H^*)=0$ & $\theta_H^*\wedge\theta_H^*=0$ \\
\hline
  $(\overline{\partial}_E)^2=0 $& $(\partial_E)^2=0$ & $F^H_{D_H}=-( \theta_H^*\wedge\theta + \theta\wedge\theta_H^*)$ \\
 \hline
\end{tabular}
\end{equation}\\
\end{corollary}

An algebraic vector bundle $E$ over $M$ is said to have a {\it
parabolic structure} if there is a collection of algebraic bundles
$E_{\alpha}$ extending $E$ over $\overline{M}$ such that the
extensions form a decreasing left continuous filtration and
$E_{\alpha+1}=E_{\alpha}\otimes \sO_M(-D_\infty).$ It is
sufficient to consider the index $0\leq \alpha <1.$ As the set of
values where the filtration jumps is discrete, it is a finite set.
Over a punctured curve $C_0=C-S,$ the {\it parabolic degree} of
$E$ is then defined by
$$\pardeg(E):=\deg \overline{E}+\sum_{s\in S} \sum_{0\leq \alpha
<1}\alpha \dim(\Gr_\alpha{\overline{E}(s)})$$ where
$\overline{E}:=E_0=\cup E_{\alpha};$ the {\it parabolic degree} of
a filtered bundle $\{E_{\alpha}\}$ on a higher dimensional $M$ is
defined just by taking the {\it parabolic degree} of the
restriction $\{ (E_{\alpha})|_C \} $ over a general curve $C_0$ in
$M$ (see the choice of $C_0$ in the following theorem \ref{DSUY}).
As any subsheaf of {\it parabolic} vector bundle $E$ has a {\it
parabolic} structure induced from $E,$ one then has the definition
of the {\it stability} for {\it parabolic} bundles (cf. \cite{Si2}). \\

A {\it harmonic bundle} $(V, H, \nabla)$ is called {\it tame} if
the metric $H$ has at most polynomial growth near the infinity.
The {\it tameness} of $(V, H, \nabla)$ is equivalent to that all
eigenvalues of the Higgs field of the induced Higgs bundle $(E,
\theta)$ have poles of order at most one at the infinity
$D_\infty.$ In other words, a {\it harmonic bundle} $(V,H,\nabla)$
is {\it tame} if and only if $(V, \nabla)$ has only regular
singularity at $D_\infty.$ Hence, any {\it tame harmonic bundle}
and its induced Higgs bundle over $M$ are algebraic. Simpson and
Jost-Zuo proved that any $\C$-local system on $M$ has a {\it tame
harmonic metric}. Moreover, if $(E,\theta)$ is a Higgs bundle
induced from a {\it tame harmonic bundle}, $E$ has a {\it
parabolic structure} compatible with the extensive Higgs field,
i.e., that one has a filtered regular Higgs bundle
$\{(E,\theta)_{\alpha} \}$ on $\overline{M}.$ Consider
$(E,\theta)$ over $\Delta^*$ as an example. Let $E_\alpha$ be an
extensive vector bundle generated by all sections $e \in
E|_{\Delta^*}$ with $|e(t)|_H\leq C|t|^{\alpha+\varepsilon} \mbox{
for } \forall \varepsilon>0.$ Then, $\theta_\alpha : E_{\alpha}\to
E_{\alpha}\otimes \Omega^1_{\Delta}(\log \ 0).$ Denote
$\overline{E}=\cup E_{\alpha}$ and $\overline{\theta}=\cup
\theta_{\alpha}.$ If all monodromies are quasi-unipotent, the
extension $\overline{E}=E_{0}$ can be chosen to be the Deligne
quasi-unipotent extension.\\

If $M$ is a punctured curve, the main theorem of Simpson in
\cite{Si2} shows that a filtered Higgs bundle $\{(E,\theta)_\alpha
\}$ is ploy-stable of {\it parabolic degree} zero if and only if
it corresponds to a ploy-stable local system of {\it degree} zero.
It can be generalized to higher dimensional bases (cf.
\cite{JZ1},\cite{Zuo1}):
\begin{theorem}[The generalized Donaldson-Simpson-Uhlenbeck-Yau correspondence]\label{DSUY}
Let $M$ be a quasi-projective manifold such that it has a smooth
projective completion $\overline{M}$ and
$D_\infty=\overline{M}\setminus M$ is a normal crossing divisor.
Let $(V,H,\nabla)$ be a {\it tame harmonic bundle} on $M$ and $\{
(E,\theta)_{\alpha} \}$ be the induced filtered Higgs bundle.
Then, one has:
\begin{myenumi}
\item $(V,H,\nabla)$ is a direct sum of irreducible ones and $\{
(E,\theta)_{\alpha} \}$ is a poly-stable filtered Higgs bundle of
{\it parabolic degree} zero.
\item If $(V,H,\nabla)$ is
irreducible,  $\{ (E,\theta)_{\alpha} \}$ is a stable filtered
Higgs bundle of {\it parabolic degree} zero.
\end{myenumi}
\end{theorem}

{\bf Remark.} Given a stable filtered Higgs bundle $\{
(E,\theta)_{\alpha} \}$ of {\it parabolic degree} zero, it is
still an open question whether there exists an irreducible {\it
tame harmonic bundle} such that $\{ (E,\theta)_{\alpha} \}$ is
induced from it. The difficulty is the existence of {\it
Hermitian-Yang-Mills metrics} on vector bundles over
quasi-projective manifolds.
\begin{proof}[Sketch of the proof]
In the higher dimensional algebraic manifold $M,$ one can choose a
smooth punctured curve $C_0$ such that its smooth completion $C
\subset \overline{M}$ is a complete intersection of very ample
divisors and it intersects $D_{\infty}$ transversally. The
homomorphism of fundamental groups $\pi_1(C_0) \onto \pi_1(M) \to
0$ is then surjective by the quasi-projective version of the {\it
Lefschetz hyperplane theorem} (cf. \cite{GM}). The restriction
$(V,H, \nabla)|_{C_0}$ is a {\it tame harmonic bundle} as Jost-Zuo
showed that the metric $H$ could be chosen to be {\it
pluriharmonic}, i.e., the metric on the restricted bundle over any
subvariety of $M$ is always {\it harmonic} \cite{JZ1,JZ2}). The
statement follows directly from Simpson's results on noncompact
curve and the subjectivity of $\pi_1(C_0) \onto \pi_1(M) \to 0.$
\end{proof}

\subsection*{Stability of Higgs bundle}
As an application, we have
\begin{theorem}\label{flat}
Let $M$ be a quasi-projective manifold such that it has a smooth
projective completion $\overline{M}$ and $D_{\infty}=\overline{M}
\setminus M$ is a simply normal crossing divisor in
$\overline{M}.$ Let $(V, H, \nabla)$ be a {\it tame harmonic
bundle} on $M$ and $(E, \overline{\partial}_E, \theta)$ be the
induced Higgs bundle. Then, $e$ is a nontrivial holomorphic
section of $(E,\overline{\partial}_E)$ with $\theta(e)=0$ if and
only if $e$ is a nonzero flat section of $(V,\nabla).$
\end{theorem}
%
\begin{proof} The \textquotedblleft$\Rightarrow$" part is obvious. We only show the
\textquotedblleft$\Leftarrow$" part.

Let $e$ be a nontrivial holomorphic section of
$(E,\overline{\partial}_E)$ with $\theta(e)=0,$ it corresponds to
the nonzero sheaf morphism $\sO_M \rightarrow \sO(E).$  Let $F$ be
the saturation sheaf generated by $e,$ $F$ is a holomorphic
subbundle of $E.$ Actually, $0\rightarrow (F,0)\rightarrow
(E,\theta)$ is a Higgs subsheaf.

{\it Step} 1. If $M$ is a curve, then $\overline{F}=\sO_M(D)$ with
$D\geq 0$ where $\overline{F}$ is the extension of $F$ to
$\overline{M}.$  Thus, $\pardeg(F)\geq 0.$ On the other hand,
Simpson showed in \cite[Lemma 6.2]{Si2} that
$$\pardeg(F)=\int_{M}\Tr(\Theta(F,H_F))$$
where $H$ is the {\it Hermitian-Yang-Mills metric} on $E$ and $H_F$
is the restricted metric on $F.$ $H$ induces a $C^{\infty}$
splitting $E=F\oplus F^{\bot},$ then
$$\Theta(F,H_F)=\Theta(E,H)|_F+\overline{A}\wedge A=-(\theta\wedge \theta^*_H)|_F
-(\theta^*_H\wedge\theta)|_F +\overline{A}\wedge A$$ where $A\in
A^{1,0}(\Hom(F,F^{\bot}))$ is the second fundamental form of
subbundle $F\subset E.$  As $\theta(F)=0,$
$\Theta(F,H_F)=-(\theta\wedge \theta^*_H)|_F+\overline{A}\wedge
A.$ Hence,
$$\pardeg(F)=\int_{M}\Tr(\Theta(F,H_F))\leq 0,$$
It is obvious that
$$\pardeg(F)=0 \Longleftrightarrow A=0 \ \mathrm{and} \ \theta^*_H|_F=0. $$
Therefore, the splitting  $E=F\oplus F^{\bot}$ is holomorphic. In
fact there is a splitting of Higgs bundle
$$(E,\theta)=(E_1,\theta)\oplus (F,0),$$ and two filtered sub Higgs
bundles are both polystable.  The {\it tame harmonic bundle} $\U$
corresponds to the Higgs bundle $(F,0),$ so it is unitary. Thus,
the metric connection on $\U$ is flat and $e$ is a nonzero flat
section of $\U$ (cf. \cite{NS} if $M$ is a compact curve).

{\it Step} 2. If $M$ is higher dimensional, we can take a generic
projective curve $C\subset \overline{M}$  such that $C$ is a
complete intersection of very ample divisors intersecting
$D_{\infty}$ transversally and $\pi_1(C\cap M) \onto \pi_1(M)$ is
surjective. Let $C_0= C\cap M,$ the restriction $H|_{C_0}$ is also
a {\it tame harmonic} metric. As in Step 1, we have a Higgs
splitting over $C_0$
$$(E|_{C_0}, \theta_{C_0})=(E_1, \theta_{C_0})\oplus (F|_{C_0}, 0)$$
where $(E_1, \theta_{C_0})$ and $(F|_{C_0}, 0)$ are Higgs subbundles
(one should be careful that the restricted Higgs field
$\theta|_{C_0}$ is not $\theta_{C_0}$). The Higgs splitting
corresponds to a local system splitting $\V=\W\oplus\U$ over $C_0$
where $\U$ is unitary part corresponding to $(F|_{C_0}, 0).$ On the
other hand, by the subjectivity of $\pi_1(C_0) \onto \pi_1(M) $ we
have a splitting of a harmonic bundle over $M$
$$\widetilde{\V}=\widetilde{\W}\oplus\widetilde{\U}$$ where $\widetilde{\V}|_{C_0}=\V,$
$\widetilde{\W}|_{C_0}=\W$ and $\widetilde{\U}|_{C_0}=\U.$
$\widetilde{\U}$ corresponds to the Higgs subbundle $(F,0)$ by the
{\it generalized DSUY correspondence}, so that it is unitary.
Then, $e$ is a flat section of $\widetilde{\U}. $
\end{proof}

{\bf Remark.} In fact, we generalize a result of Simpson on
compact manifolds (cf. \cite[Lemma 1.2]{Si3}). The result can also
follow from the Bochner method (cf. \cite{Sch}) by using the
non-positivity of curvature forms and the estimates of curvatures
at the infinity (cf. \cite{JZ1,JZ2}).

\section{The Geometry of Lefschetz Pencils}
\subsection*{Lefschetz pencils}
\begin{definition}[cf. \cite{Del1}]\label{lefschetz}
Let $X\subset \P^N$  be a projective manifold and  $L$ be a fixed
hyperplane in $\P^N$ of dimension $N-2.$ The set of hyperplanes $
H_s \subset \P^N$ passing through $L$ is parameterized by a
projective line $\P^1$ in the dual space $(\P^N)^*\simeq \P^N.$
$L$ can be chosen to satisfy that:
\begin{enumerate}
    \item [(a)] $L$ intersects
$X$ transversally so that $Y=X\cap L $ is a smooth subvariety of
$X.$
    \item [(b)] There exists a finite subset $S=\{s_1, \cdots s_k \} \subset
\P^1$ such that
\begin{itemize}
    \item for any $s_j \in S,$ the variety $X_{s_j}$ has one and only one ordinary double
 singularity $x_j\in Y\cap X_{s_j};$
    \item if $s \notin S,$ $H_s$ intersects $X$ transversally
(then $X_s=X\cap H_s $ is a nonsingular variety).
\end{itemize}
\end{enumerate}
Such a family $\{X_s\}_{s\in \P^1}$ is called {\it Lefschetz
pencil}.
\end{definition}

{\bf Remark.} For a general projective line $C\subset (\P^N)^*,$
$L=\bigcap_{s\in C} (X\cap H_s)$ actually satisfies all conditions
in the definition.
\begin{lemma}[cf. \cite{Del1}]
Any Lefschetz pencil is a proper flat family with a section.
\end{lemma}
\begin{example}[A Typical Lefschetz Pencil of Calabi-Yau Varieties]\label{Lefschetz-CY}
Let $(Y,\sL)$ be a projective manifold with a very ample line
bundle. The Veronese embedding $\iota: Y \hookrightarrow
\P|\sL|=\P$ is given by the complete linear system
$|\sL|=H^0(Y,\sL),$ then $\sL =\iota^*\sO_{\P}(1)$ and
$H^0(\P,\sO_{\P}(1))=H^0(Y,\sL).$ Suppose that $Y=\P^{n+1}$ and
$\sL=\sO_{Y}(n+2).$ The Veronese embedding becomes $\iota :
\P^{n+1} \hookrightarrow \P^{N}$ where
$N=\left(\begin{array}{ccccccccc}
2n+3  \\
n+1  \\
\end{array}\right)-1.$ 
Blowing up $X:=\iota(\P^{n+1})$ centered at a certain codimension
$2$ subvariety, one obtains a nonsingular projective variety
$\widetilde{X}$ and a Lefschetz pencil $f: \widetilde{X}
\rightarrow \P^1 $ of $n$-dimensional Calabi-Yau varieties.
\end{example}

Let $F, H$ be smooth degree $5$ homogenous polynomials with $H$ is
not in the Jacobian ideal of $F.$ The family $t_0F + t_1H$ by the
parameter $[t_0, t_1]\in \P^1 $ may not be Lefschetz pencil, but
is it rigid? 
\begin{proposition}\label{number-singularity}
Let $f:\ssX \rightarrow \P^1$ be a family. Assume that the {\it
infinitesimal Torelli theorem} holds for a general fiber and $f$
is not a local constant analytic family. Then, $f$ has at least
$3$ singular fibers.
\end{proposition}
\begin{proof}[Sketch of the proof]
It is sufficient to show that the restricted family $f: \ssX \to
\C \setminus\{0 \} $ is not smooth. Replacing $\C\setminus \{0\}$
by an unramified covering:
$$ (\C\setminus \{0\},z )\> t=z^k >> (\C\setminus \{0\},t),$$
we have that the pull back family $f'$ has only unipotent
monodromy VHS. Then, the monodromy of the VHS must be trivial
because the global monodromy is semi-simple by Deligne's complete
reducible theorem (cf. \cite{Del}). Thus, $f$ is a local constant
analytic family.
\end{proof}
\begin{corollary}\label{CY-singularity} Let $f:\ssX \rightarrow \P^1$
be a Lefschetz pencil of Calabi-Yau varieties such that the period
map is injective at one point of $\P^1,$ i.e., non-isotrivial,
then $f$ has at least $3$ singular fibers.
\end{corollary}

{\bf Remark.} Actually, one has a general result: Consider a
nonisotrivial family $(f: \sX \rightarrow \P^1).$ If $X$ is a
projective manifold of
non-negative Kodaira dimension, then $f$ has at least $3$ singular fibres (cf. \cite{VZ1}). \\
\subsection*{Vanishing cycles spaces of Lefschetz pencils}
Let $f:\ssX \rightarrow \P^1 $ be a Lefschetz pencil of
$n$-dimensional varieties degenerate at $S=\{s_1, \cdots, s_k \}.$
Let $\Delta_i$ ($1\leq i \leq k$) be a set of disjoint disks
centered at the points $s_i $ ($1\leq i \leq k$) each other and
let $\Delta^*_i=\Delta_i \setminus \{s_i\}.$ Let $f_{1}:
\ssX_{1}=f^{-1}(\Delta_1) \rightarrow \Delta_1$ be the restricted
family, there is a unique singularity $x \in f_{1}^{-1}(0)$ which
is simple. If one chooses a holomorphic local coordinates
$z=(z_0,\cdots,z_n)$ in a suitable neighborhood of
$x=(0,\cdots,0),$ $f_1$ can be written as $$f_1(z)=z_0^2+\cdots +
z_n^2.$$ Let $B_{\epsilon}(x)$ be a small ball in $\ssX_{1}$
centered at $x.$ If $s\in\Delta$ is sufficiently closed to $0,$
the homotopy type of $V_s:=B_s\cap \ssX_s$ is same as a real
$2n$-dimensional sphere, so $H^n_c(V_s,\Z)\simeq \Z$ and there is
a natural inclusion
$$\iota : H^n_c(V_s,\Z) \rightarrow H^n(\ssX_s,\Z). $$
The image $\delta=\iota(1)$ (up to multiplying $-1$) is called
{\it vanishing cycle}. Similarly, one has $f_i:\ssX_i \to
\Delta_i$  and vanishing cycles $\delta_i$ in $H^n(\ssX_{s'_i})$
where $s'_i$ is near $s_i$ for $2\leq i\leq k.$ Fixing a point
$s_0 \in \P^1 \setminus S,$ there is a monodromy transformation
identifying $H^n(\ssX_{s'_i},\Z)$ with $H^n(\ssX_{s_0}, \Z)$ for
each $1\leq i\leq k.$ Hence, in this way $\delta_i \in
H^n(\ssX_{s_0}, \Z)$ for all $i.$

The $\Q$-liner space $V$ generated by all $\delta_i \in
H^n(\ssX_{s_0}, \Z)$ is called {\it vanishing cycles space}. One
has the Picard-Lefschetz formula (cf. \cite{Gro}):
\begin{equation}\label{Picard-Lefschetz}
 T_i(\upsilon)=\upsilon+ (-1)^{(n+1)(n+2)/2}
(\upsilon, \delta_i)\delta_i, \, \,  \upsilon\in
H^n(\ssX_{s_0},\Q)
\end{equation}
with
\begin{equation*}
(\delta_i, \delta_i)= \left\{
\begin{array}{ll}
     0 &  \hbox{ $n$ odd } \\
     2(-1)^{n/2} &  \hbox{ $n$ even.}
\end{array}
\right.
\end{equation*}
The transform $T_i :H^n(\ssX_{s_0}, \C)\to H^n(\ssX_{s_0}, \C)$ is
generated by a loop around $s_i$ and is called a {\it
Picard-Lefschetz transformation}. The nondegenerated intersection
form $(,)$ is preserved by all {\it Picard-Lefschetz
transformations}. If one considers only a local Lefschetz pencil
$f_{1}: \ssX_{1} \rightarrow \Delta_1,$ the local monodromy
transform $T:H^n(\ssX_t,\C)\to H^n(\ssX_t,\C)$ satisfies that:
\begin{equation*}
\begin{array}{ll}
     (T-\id)^2=0 &  \hbox{ $n$ odd } \\
\  \  \ \ \ \ \ \ T^2=id &  \hbox{ $n$ even. }
\end{array}
\end{equation*}
Thus, the Landman theorem holds naturally for Lefschetz pencils.
\begin{lemma}[Lefschetz in the classical case cf. \cite{Del1}]\label{class-Lefschetz} Vanishing
cycles are conjugate under the action of $\pi_1(\P^1 \setminus S,
s_0)$ up to sign. In particular, the vanishing cycles space $V$ is
stable under the action of $\pi_1(\P^1 \setminus S, s_0).$
\end{lemma}
\begin{proposition}\label{CY-vanishing-cycle}
Let $f:\ssX \rightarrow \P^1$ be a Lefschetz pencil. Assume that
the {\it infinitesimal Torelli theorem} holds for generic fiber.
If the Torelli map of $f$ is injective at one point of $\P^1,$
i.e., $f$ is non-isotrivial, the vanishing cycles space $V$ would
be non trivial.
\end{proposition}
\begin{proof}
Let $f^0:\ssX_0=f^{-1}(\P^1\setminus S) \to \P^1\setminus S$ be
the maximal smooth subfamily of $f.$ The set $S$ is non-empty and
$\# S \geq 3$ by \ref{number-singularity}.  Suppose $V=0, $ then
$T_j\equiv id \ \forall j$ by the Picard-Lefschetz formula and the
VHS $R_{prim}^nf^0_{*}\Q$ on $\P^1 \setminus S$ can be extended to
a VHS on $\P^1.$ But the Torelli map becomes a constant map
because $\P^1$ is simply connected.
\end{proof}

\subsection*{Invariant subspace of cohomology group }
Let $f:\ssX\rightarrow Y$ be a  smooth proper family over an
algebraic manifold $Y$ and $t_0\in Y$ be a fixed point. Denote
$X=f^{-1}(t_0)=\ssX_{t_0}.$ For each $t\in Y,$ $\pi_1(Y, t_0)$
acts on $H^n(X_t, \Q)$ naturally and one has a $\pi_1(Y,
t_0)$-invariant subspace $$H^n(\ssX_t, \Q)^{\pi_1(M,
t_0)}\hookrightarrow H^n(\ssX_t, \Q)$$ such that the inclusion is
a Hodge morphism of type $(0, 0).$ One then has a constant sheaf
$(R^nf_*(\Q))^{\pi_1(Y, t_0)}$ by gluing  $\{H^n(\ssX_t,
\Q)^{\pi_1(Y, t_0)}\}_{t\in Y}$ together. Moreover, for $\forall t
\in Y$ there is a natural Hodge isomorphism $$\kappa_t: H^0(Y,
R^nf_*(\Q)) \> \simeq
>> H^n(\ssX_t, \Q)^{\pi_1(Y, t_0)}.$$

\begin{theorem}[Deligne \cite{Del}]\label{deligne}
Let $f:\ssX\rightarrow Y$ be a  smooth proper family over non
singular algebraic variety $Y.$  One has:
\begin{enumerate}
    \item Leray's spectral sequence for $f$
$$E^{p,q}_2=H^p(Y,R^qf_*(\Q))\Rightarrow H^{p+q}(\ssX,\Q)$$
degenerates at $E_2.$
    \item Let $\overline{\ssX}$ be any smooth
compactification of $\ssX$ and $i:\ssX \hookrightarrow
\overline{\ssX}$ be the inclusion. Then, the composite Hodge
morphism
$$H^n(\overline{\ssX},\Q)\rightarrow H^n(\ssX,\Q) \rightarrow H^0(Y,R^nf_*(\Q))$$ is
surjective.
    \item $(R^nf_*(\Q))^{\pi_1(Y,t_0)} $ is the maximal constant sub $\Q$-VHS of $R^nf_*(\Q).$
\end{enumerate}
\end{theorem}

The theorem \ref{deligne} says that the $(p,q)$-component
$\omega^{p,q}$ of a global section $\omega$ of a VHS is invariant
under $\pi_1(Y,t_0),$ so $\omega^{p,q}$ is also a global section.
Therefore, if a global section $\omega$ of
$(R^nf_*(\C))^{\pi_1(Y)}$ is of $(p,q)$-type at one point of a
connected manifold $Y,$ it is of $(p,q)$-type everywhere.
\begin{corollary}
Let $f: \ssX \rightarrow \P^1$ be a Lefschetz pencil of $n$-folds
smooth over $C_0=\P^1 \setminus S.$ Let $s_0 \in C_0$ be a fixed
point and $V$ be the vanishing circle space of
$H^n(\ssX_{s_0},\Q)$ and $V^{\bot}$ be the orthogonal complement
of $V$ in $H^n(\ssX_{s_0},\Q)$ under the nondegenerated
intersection $(,).$ Assume that $V$ is not trivial, then
$V^{\bot}=H^n(\ssX_{s_0},\Q)^{\pi_1(C_0, s_0)}.$
\end{corollary}
The proof depends on the Picard-Lefschetz formula. Altogether,
\begin{theorem}[Deligne \cite{Del1}]\label{irrducible-module} Let
$f:\ssX\rightarrow \P^1$ be a Lefschetz pencil of $n$-dimensional
varieties and $f^0:\ssX_0\rightarrow \P^1 \setminus S$ be the
maximal smooth subfamily of $f.$  Let $\V$ be a $Z$-local system
generated by the stable action of $\pi_1(\P^1\setminus S, s_0)$ on
the vanishing cycles space $V$ where $s_0\in \P^1\setminus S$ is a
fixed base point. Then, one has a $\Q$-VHS splitting
$$ R^nf^0_{*}(\Q)=\V_\Q \oplus (R^nf^0_{*}(\Q))^{\pi_1(\P^1\setminus S, s_0)} .$$
Moreover, if $\V$ is nontrivial then it is absolutely irreducible,
i.e., for any algebraic field extension $\K /\Q,$ $\V_{\K}$ is an
irreducible $\pi_1(\P^1\setminus S,s_0)$-module.
\end{theorem}

{\bf Remark.} The complete irreducible theorem of Deligne induces
a $\Q$-local system decomposition $R^nf^0_{*}\Q=\oplus_i \V_{i
\Q}.$  Each $\V_{i\Q}$ has an induced Hodge filtration from the
Hodge filtration of $R^nf^0_{*}\Q.$ Hence, the decomposition is
actually of VHS, but in
general it is not compatible with polarizations. \\

\subsection*{Lefschetz pencils of Calabi-Yau Varieties}
\begin{lemma}\label{zhangyi1}
Let $f:\ssX \rightarrow Y$ be a smooth family of Calabi-Yau
$n$-folds and $\sI$ be the Higgs bundle induced from the $\Q$-VHS
$(R^nf_*(\Q))^{\pi_1(Y)}.$ If the differential of the Torelli map
of $f$ is injective at some points in $Y,$ then
$$\dim_{\C}\sI^{n,0}=\dim_{\C}\sI^{0,n}=0.$$
\end{lemma}
\begin{proof}
Suppose that the differential of the Torelli map of $f$ is
injective at $y_0\in Y.$ Let $\phi: (\ssX, X) \rightarrow (S, 0)$
be the Kuranishi family of $X=f^{-1}(y_0).$ The theorem of
Bogomolov-Todorov-Tian says that the deformation of $X$ is
unobstructed, so $S$ is smooth (cf. \cite{Ti},\cite{To89}). We
then have a commutative diagram over a neighborhood $U$ (in the
topology of complex analytic spaces) of $y_0$ in $Y$
$$
\begin{TriCDV}
{(U, t_0)} {\> \pi >>}{(S, 0) } {\SE \lambda_U EE}{\SW W\lambda W}
{\prod \Gr(h^p, H_\C)}
\end{TriCDV}
$$
where $\lambda, \lambda_U$ are Torelli maps and $\pi$ is the
Kuranishi map. Moreover, we can assume that $\lambda_U$ is an
embedding by contracting $U$ sufficiently.

We claim that
$$h^{n, 0}(H^n(\ssX_{t_0}, \C)^{\pi_1(Y, t_0)})<
h^{n,0}(\ssX_{t_0})=1 .
$$
Otherwise,
$h^{n,0}(H^n(\ssX_{t_0},\C)^{\pi_1(Y,t_0)})=h^{n,0}(\ssX_{t_0})=1.
$ Let $\overline{\ssX}$ be the smooth Haronaka compactification of
$\ssX,$ there are composite morphisms:
\begin{multline*}
\overline{i}^*_{t}: H^n(\overline{\ssX},\C)\> j\circ i^*
  >> H^0(Y,R^nf_*\C) \> \kappa_{t}
>> H^n(\ssX_{t},\C)^{\pi_1(M,t)} \\
\hookrightarrow H^n(\ssX_{t},\C) \, \forall t \in Y,
\end{multline*}
and each $\overline{i}^*_{t}$ is a morphism of Hodge structure of
type $(0,0)$ (cf. \cite{Del}). By \ref{deligne},
$$h^{n,0}(H^n(\ssX_{t},\C)^{\pi_1(Y,t)})= h^{n,0}(\ssX_t)=1 \, \forall
t \in Y. $$ Therefore, each $\overline{i}^*_{t}$ has to be
surjective and all $(n,0)$ holomorphic forms of $H^n(\ssX_{t},\C)$
lift to $H^n(\overline{\ssX},\C),$ i.e., that we have a  Hodge
isomorphism
$$ H^n(\ssX_t,\Q) \> \simeq >>  H^n(\ssX_t,\Q)^{\pi_1(Y,t_0)} \> \simeq >>  H^0(Y,R^nf_*(\Q)).$$
It is a contradiction to that $\lambda_U$ is an embedding over
$U.$
\end{proof}
\begin{corollary}\label{WP-VHS}
Let $f:\ssX \rightarrow Y$ be a proper smooth family of Calabi-Yau
$n$-folds.  If the differential of the Torelli map is injective at
some points, then $f_*(\omega_{\ssX/Y})$ and its multi-tensors are
geometric nontrivial.
\end{corollary}
\begin{proof}
Suppose that the differential of the Torelli map of $f$ is
injective at $y_0\in Y.$ Let $\pi:\sX \rightarrow \ssM $ be the
maximal subfamily of the Kuranishi family for $X=f^{-1}(y_0)$ with
respect to a fixing polarization. The polarized Kuranishi base
$\ssM$ is smooth and the Kodaira-Spencer map $\rho(t):T_{\ssM,
t}\to H^1(\sX_t,T_{\sX_t})$ is injective  for $\forall t \in
\ssM.$ In the category of complex analytic spaces, the polarized
Kuranishi family is universal. Actually, under the assumptions we
have a commutative diagram
$$
\begin{CDS}
\ssX_U \> \hookrightarrow >> \sX  \\
\V V f_U V \novarr \V V \pi V  \\
U \> \hookrightarrow >> \ssM
\end{CDS}
$$
where  $f_U:\ssX_U= \ssX\times_{Y}U \to U$ is the restricted
family over a small open neighborhood  $U\subset Y$ of $y_0.$
Hence, we have
$$\pi_*(\omega_{\sX/\ssM})|_{U}=\pi_*(\omega_{\sX/\ssM}|_{\ssX_U})=f_*\omega_{\ssX/Y}|_{U}.$$
On the other hand, on $\ssM$
$$\omega_{WP}=-\frac{\sqrt{-1}}{2}\partial\overline{\partial}\log
 h =c_1(\pi_*\omega_{\sX/\ssM}, h)$$  where $h$
is the Hodge metric and $\omega_{WP}$ is the K\"ahler form of the
Weil-Petersson metric (cf. \cite{Ti},\cite{To89}). Therefore,
$$\int_{U}(c_1(f_*\omega_{\ssX/Y},h))^{\dim
Y}=\int_U (c_1(\pi_*(\omega_{\sX/\ssM}),h))^{\dim Y}>0.$$
\end{proof}

A polarized VHS over $Y$ is {\it isotrivial} if it becomes a
constant VHS after a finite \'etale base change.  It is obvious
that a polarized VHS is {\it isotrivial} if and only if the Hodge
filtration is local constant (cf. \cite{Del}).

\begin{proposition}\label{zhangyi2}
Let $f:\ssX\rightarrow \P^1$ be a Lefschetz pencil of $n$-folds
smooth over $P^1\setminus S$ and $f^0:\ssX_0\rightarrow \P^1
\setminus S$ be the maximal smooth subfamily of $f.$  Assume that
$n$ is even, then the VHS $R^nf^0_{*}(\Q)$ must be isotrivial and
$$(R^nf^0_{*}(\C))^{n,0},(R^nf^0_{*}(\C))^{0,n} \subset
(R^nf^0_{*}(\C))^{\pi_1(\P^1\setminus S)}. $$
\end{proposition}
\begin{proof}
All Picard-Lefschetz transforms are of order $2$ as $n$ is even.
Since $\P^1$ is topologically simply-connected, the global
monodromy group $\Gamma$ is a finite commutative group. Actually,
$\Gamma \cong (\prod_1^{\#S} \Z_2,+).$ The VHS $R^nf^0_{*}(\Q)$
must be {\it isotrivial} by \cite[Theorem 9.8]{Gr2}. Then,
$$T_i (H^{n,0}(\ssX_{s_0},\C)) \subset H^{n,0}(\ssX_{s_0},\C) \  \forall i
$$ where $s_0$ is a fixed base point in $\P^1\setminus S.$
By the Picard-Lefschetz formula, we have $(\eta,\delta_i)\delta_i
\in H^{n,0}(\ssX_{s_0},\C)\ \forall i$ for any $\eta \in
H^{n,0}(\ssX_{s_0},\C).$

We claim that $(\eta,\delta_ j)=0 \ \forall j $ for any $\eta \in
H^{n,0}(\ssX_{s_0},\C).$ Otherwise, $(\eta,\delta_ j)\neq 0$ for
one $j$ and one $\eta\neq 0,$ then $\delta_j\in
H^{n,0}(\ssX_{s_0},\C).$ But it is a contradiction to that
$(\delta_j, \delta_j)\neq 0.$ From the formula
\ref{Picard-Lefschetz}, we then have:
$$H^{n,0}(\ssX_{s_0},\C) \subset
H^n(\ssX_{s_0},\C)^{\pi_1(\P^1\setminus S, s_0)}.$$
\end{proof}
\begin{corollary}
Let  $f:\ssX \to \P^1$ be a Lefschetz pencil  of Calabi-Yau
varieties. Assume that the Torelli map is injective at some
points, the global algebraic monodromy group is infinite.
\end{corollary}

Altogether, from \ref{irrducible-module}, \ref{zhangyi1} and
\ref{zhangyi2} we obtain a key result for the {\it rigidity} of
Lefschetz pencils of Calabi-Yau varieties.
\begin{theorem}\label{zhangyi3}
Let $f:\ssX \to \P^1$ be a Lefschetz pencil of $n$-dimensional
Calabi-Yau varieties and $f^0:\ssX_0=f^{-1}(\P^1 \setminus S )\to
\P^1 \setminus S$ be its maximal smooth subfamily. Let $\V$ be the
$\Z$-local system generated by the vanishing cycles space. If the
differential of the Torelli map of $f^0$ is injective at some
points (for example, non-isotrivial families). Then, we have:
\begin{myenumi}
\item The integer $n$ must be odd.

\item $\V$ is an nontrivial absolutely irreducible $\pi_1(\P^1
\setminus S)$-module and
$$(R_{prim}^nf^0_{*}(\C))^{n,0}, (R_{prim}^nf^0_{*}(\C))^{0,n}
\subset \V_\C.$$

\item  $\V_\Q$ is the $\Q$-sub local system of
$R_{prim}^nf^0_{*}(\Q).$ Moreover,
$$R_{prim}^nf^0_{*}(\Q)=\V_\Q \oplus (R_{prim}^nf^0_{*}(\Q))^{\pi_1(\P^1\setminus S)}$$
by the relative Lefschetz decomposition.
\end{myenumi}
\end{theorem}

\section{A Criterion for Rigidity and Its Applications} 
\subsection*{Endomorphisms of Higgs bundles over a product
variety}\label{product-in-the-section} Let $S_0, T_0$ be
quasi-projective manifolds such that they have smooth projective
completions $S,T$ and $D_S=S-S_0, D_T=T-T_0$ are normal crossing
divisors. Let $(\V,\nabla)$ be an arbitrary polarized $\R$-VHS
over $S^0\times T^0$ such that local monodromies around  the
divisor $D$ at infinity are quasi unipotent. Let $(E,\theta)$ be
the Higgs bundle induced from  $\V.$ Extending the Higgs bundle to
the infinity, we have the quasi canonical extension
$(\overline{E},\overline{\theta})$ with
$$\theta: \overline{E} \rightarrow \overline{E}\otimes \Omega^1_{S \times T}(\log D).$$
We then obtain a flat endomorphism of $\mathbb{V}|_{\text{S}_{t}}$
with Hodge type $(-1,1)$ (cf. \cite{JY} and \cite{Zuo0}) in the
following way.
\begin{itemize}
    \item First, one has a decomposition
$$ \Omega_{S\times T}^{1}(\log D)=p_{S}^{*}\Omega _{S}^{1}(\log
D_{S})\oplus p_{T}^{*}\Omega_{T}^{1}(\log D_{T})$$ where $ p_{S}:
S\times T \rightarrow  S$ and $p_{T}: S \times T \rightarrow T$ are
projections. The Higgs map
\begin{equation}\label{differential-of-period-map}
\theta : p_{S}^{*}\Theta_{S}(-\log D_{S})\oplus
p_{T}^{*}\Theta_{T}(-\log D_{T})\rightarrow End(\overline{E})
\end{equation}
is just the sheaf map of the differential of the period map. Fix a
$t \in T^0,$ the restricted Higgs map over $S_{t}$ is
\begin{equation}
\theta|_{S_{t}}:(p_{S}^{*}\Theta_{S}(-\log D_{S})\oplus
p_{T}^{*}\Theta_{T}(-\log D_{T}))|_{S_{t}}\rightarrow End(\overline{E})|_{S_t}. \label{kang}%
\end{equation}
Note that
$$ p_{T}^{*}\Theta_{T}(-\log D_{T})|_{S_t} \simeq \bigoplus^l
\mathcal{O}_{S_t}$$ where $l$ is the dimension of $T.$ Let $
1_{T}$ be a constant
section of $p_{T}^{*}\Theta_{T}(-\log D_{T})|_{S_t}.$ Then, we obtain an endomorphism%
\begin{equation}
\sigma:=\theta|_{S_t}(1_{T}):E|_{S_t^{0}} \rightarrow E|_{S_t^0}.
\label{JY-Kang}
\end{equation}
and $\sigma$ must be of $(-1,1)$-type. Moreover, $\sigma$ is a
morphism of Higgs sheaf, i.e., there is a commutative diagram
$$
\begin{CD}
E|_{S_t^0} \> {\theta_{S^0_t}}>> E|_{S_t^0} \otimes\Omega^1_{S^0_t} \\
\V V \sigma V           \V V {\sigma\otimes id} V \\
E|_{S_t^0} \> {\theta_{S^0_t} }>> E|_{S_t^0}\otimes\Omega^1_{S^0_t}
\end{CD}
$$ by
$\theta_{S^0_t}\wedge\theta_{S^0_t}=0,$ where $$ \theta_{S_{t}}:
\Theta_{S_t}(-\log D_{S_t})\to (p_{S}^{*}\Theta_{S}(-\log
D_{S})\oplus p_{T}^{*}\Theta_{T}(-\log D_{T}))|_{S_{t}}\to
\End(\overline{E})|_{S_t}
$$ is the Higgs field of $E|_{S_t}.$

    \item Let $M$ be any quasi-projective manifold with a smooth compactification $\overline{M}$ and
    a normal crossing divisor $D_{\infty}=\overline{M}-M. $
    Let $(E,\theta)$ be a Higgs bundle on $M.$ If $E$ carries $\R$-structure,
    same as \cite[Lemma 2.11]{Si3} for compact $M,$ the dual Higgs bundle is
$$(E^{\vee}=\bigoplus_{p+q=n}E^{\vee -p, -q}, \theta_{\vee})$$
where $E^{\vee -p, -q}=(E^{p, q})^{\vee}=E^{q, p}, $
$\theta_{\vee}^{-p,-q}=-\theta^{q, p}$ and
$$\theta_{\vee}^{-p, -q}: E^{\vee -p, -q}\rightarrow E^{\vee -p-1, -q+1}\otimes \Omega_M^1.$$
(However, for a polarized $\R$-VHS $\V_\R$ of weight $n$ one has a
natural nondegenerated
pairing $\V_\R \times \V_\R \rightarrow \R(-n)$, hence $\V_\R^\vee=\V_\R(n)$). 
Therefore, $\End(E)$ is a Higgs bundle, i.e.,
$$(\End(E), \theta^{end})=
(\bigoplus_{(p-p')+(q-q')=0}E^{p, q}\otimes E^{\vee -p',-q'},
\theta^{end})$$ with $\theta^{end}(u\otimes
v^{\vee})=\theta(v)\otimes v^{\vee}+u\otimes
\theta_{\vee}(v^{\vee}).$ All local monodromies for $\End(E)$ are
quasi unipotent as so are for $E,$ and one has Deligne's
quasi-canonical extension for $\End(E).$ Moreover, one has:
\begin{lemma}[Proposition 2.1 in \cite{Zuo0}]\label{kang2}
Let $(\End(E), \theta^{end})$ be the Hodge bundle corresponding to
a polarized VHS on $\End(\V_\R)$ which is induced by the polarized
$\R$-VHS $\V_\R.$ Then,
$$\theta^{end}(d\phi(T_{\overline{M}})(-\log D_\infty))=0.$$
\end{lemma}
    \item  The image of the map
\ref{kang} is  a trivial Higgs subsheaf as it is contained in the
kernel of a Higgs field, and $\sigma$ is a flat section by the
poly-stability of the Higgs sheaf. Precisely,
\begin{proposition}\label{sigma-flat} The endomorphism $\sigma$
obtained above is a flat $(-1,1)$-type section of
$\End(E)|_{S^0_t}. $ Therefore, $\sigma$ is an endomorphism of a
$\C$-local system, i.e.,
\begin{equation}
\sigma: \V_{\C}|_{S_t^{0}} \rightarrow \V_\C|_{S_t^0}.
\label{JY-Kang2}
\end{equation}
\end{proposition}

\begin{proof}  

First, $\theta|_{S^0\times T^0}=d\phi $ where $\phi : S^0\times
T^0 \to  D/\Gamma$ is the period map associated to $\V.$ Denote
$\Image(\theta)$ is the image sheaf. We have
$$\sigma_{s,t}\in (d\phi(\Theta_{S^0\times
T^0}))_{s,t}=(\Image(\theta))_{s,t}$$ and
$\theta^{end}_{s,t}(\sigma_{s,t})=0$ for all $(s,t)\in S^0_t$ by
\ref{kang2}.

On the other hand, we can else get $\sigma$ if we fix $t$ and vary
$s\in S^0.$ Thus, $\sigma$ is of type $(-1,1)$ and
$\theta^{end}_{S^0_t}(\sigma)=0.$ By \ref{flat}, $\sigma$ is a
flat section of $(\End(\V_\C)|_{S^0_t},\nabla).$
\end{proof}
\end{itemize}
\subsection*{ A criterion for infinitesimal rigidity}\label{infi-rigid}
In the following sections we study the manifolds for which the
{\it infinitesimal Torelli theorem} holds, and we assume that
there is a natural condition for a smooth family $f:\ssX \to Y$ of
projective $n$-folds:
\begin{center}
($\ast \ast$) \ \  {\it The differential of the period map
   for the VHS $R_{prim}^nf_*(\Q)$ is injective at some points of $Y.$}
\end{center}
By the {\it infinitesimal Torelli theorem}, the condition is
equivalent to that the induced moduli map $\eta_f: Y \to
\mathfrak{M}_h $ is a generic finite morphism, i.e., $f$ contains
no isotrivial subfamily such that the base is a subvariety passing
through a general point of $Y.$ Suppose that $Y$ is a curve, the
($\ast \ast$)-condition is then equivalent to that $f$ is a
non-isotrivial family.

\vspace{0.2cm}

A smooth family $f$ is {\it rigid} if there exists no non-trivial
{\it deformation} of $f$ over a non-singular quasi-projective
curve. A {\it deformation} of a smooth family $f$ over a
quasi-projective variety $T^0$ with  base point $0\in T^0$ is a
smooth projective morphism $g:\sX \to Y\times T^0$ with a
commutative diagram
$$
\begin{CD}
\ssX \> \simeq >> g^{-1}(Y\times\{0\}) \>\subset >> \sX \\
\V f VV \V V V \V V g V \\
Y \> \simeq >> Y\times\{0\} \> \subset >> Y\times T^0
\end{CD}.
$$
A family of varieties is {\it rigid} if its maximal smooth
subfamily is {\it rigid}.

In the category of complex analytic spaces, we replace $T^0$ by a
small disk. Then, we should say {\it infinitesimal rigidity}
instead of {\it rigidity}. It is obvious that an algebraic family
is automatically {\it rigid} if it is {\it infinitesimal rigid} in
the category of complex analytic spaces. If there is no ambiguous,
we often do not distinguish between the two notations.

Let $f:\ssX \to Y$ be a smooth polarized family of projective
$n$-folds. If a deformation $g$ of $f$ over a quasi-projective
smooth curve $T^0$ is not trivial, then the period map for $g$ is
not degenerate along $T^0$-direction at some points of $Y\times
\{0\}$ by the {\it infinitesimal Torelli theorem}. Therefore, with
our methods in studying Higgs bundles over a product variety we
reprove a theorem of Faltings, Peters and Jost-Yau (cf.
\cite{Falt},\cite{Pe},\cite{JY}).
\begin{theorem}[A criterion for rigidity]\label{rigidity}
Let $f:\ssX \to Y$ be a smooth family of polarized projective
$n$-folds satisfying the $(\ast \ast)$ condition. If $f$ is
nonrigid, there exists a flat nonzero section $\sigma$ of
$$\End(R_{prim}^nf_*(\C))^{-1,1}.$$
Moreover, the Zariski tangent space at $[f]$ of the deformation
space of $f$ is into $\End(R_{prim}^nf_*(\C))^{-1,1}. $
\end{theorem}

Saito-Zucker and Zuo also obtained similar criterions (cf.
\cite{SZ},\cite{Zuo0}). It is not difficult for us to generalize
the criterion to nonrigid polarized VHSs with this method.

Notably the non-positivity of curvatures in the horizontal
directions of period domains will underly the validity of this
criterion of the rigidity. Here, we explain the role played by the
differential geometry of period domains: Consider the simplest
case. Let $f_t:\ssX_t \to Y$ be one parameter nontrivial
holomorphic deformations of  families of polarized Calabi-Yau
manifolds satisfying the $(\ast\ast)$-condition. We then have one
parameter period maps $g_t$ from $Y$ into a fixed period domain.
The infinitesimal deformation of $g_t$ gives rise to a holomorphic
section of a Hermitian holomorphic vector bundle with non-positive
curvature in the sense of Griffiths. By the Bochner method and the
estimates of curvatures at the infinity, this holomorphic section
must be parallel.

It is the original ideal of Jost-Yau to deal with the {\it
rigidity} of Shafarevich problems (cf. \cite{JY}). We shall point
out that the $(-1,1)$-type of this flat section is a key point in
proving the {\it rigidity} throughout this paper, and it seems one
can not obtain these results only by the pure techniques of the
differential geometry. It is successful for us to use the theory
of algebraic Higgs bundles.

\subsection*{Applications of the criterion}
Let $f:(\ssX,\sL) \to Y$ be a smooth family of polarized
Calabi-Yau $n$-folds with $(X,L)=f^{-1}(0)$ and let $\pi:\ssY
\rightarrow \ssM_{c_1(L)} $ be the maximal subfamily of Kuranishi
family of $\pi^{-1}(0)=X$ with a fixed polarization $L.$ For each
$t \in \ssM_{c_1(L)},$ the Kodaira-Spencer map
$$\rho(t):T_{\ssM,t} \> \simeq >>
H^1(\ssX_t,T_{\ssX_t})_{c_1(L)}$$ is isomorphic. Let $\Phi$ be the
period map of $\pi$ and $U$ be a sufficient small neighborhood (in
topology of complex analytic spaces) of $0$ in $M.$ If the
differential of the period map $\Psi : Y \>> > D/\Gamma$ is
injective at $0\in Y,$ one has a locally commutative
diagram  
$$
\begin{TriCDV}
{(U,0)} {\> \pi >>}{(\ssM_{c_1(L)},0)} {\SE \Psi \ EE}{\SW W \Phi
W} {(D,0)}
\end{TriCDV}
$$
where $\pi$ is an embedding over $U.$ Let $(E,\theta)$ be the
Higgs bundle induced from $R_{prim}^nf_*(\Q),$ the horizonal
tangent space $T^h_{D,0}$ is contained in
$$\bigoplus_{p=0}^{n}\Hom(E^{n-p,p}|_{0},E^{n-p-1,p+1}|_{0})$$ and
$d\Psi|_0=(\theta^{n,0}|_0,\cdots,\theta^{0,n}|_0)$ where
$$\theta^{p,q}: E^{p,q} \to E^{p-1,q-1}\otimes \Omega^1_Y  \mbox{ with } p+q=n $$
is the $(p,q)$-component of the Higgs field  $\theta.$
\begin{proposition}[The infinitesimal Torelli theorem]\label{CY-sigma}
The differential of the period map of a family of Calabi-Yau
$n$-folds is injective at one point if and only if $\theta^{n, 0}:
E^{n,0} \to E^{n-1,1}$ is injective at this point.
\end{proposition}

If $f$ is nonrigid, one would has a non-trivial deformation $g$
over $Y\times T$ such that $g$ is not degenerate along the
orientation of $T$ at some points $\{s_0,\cdots,s_k \}$ of
$(Y,0)\subset Y\times T^0.$ We then have a nonzero $(-1,1)$-type
endomorphism $\sigma$ of $E$ which is flat under the Gauss-Manin
connection. As $\sigma$ is obtained from the period map over
$Y\times T^0$ along the $T^0$-direction, $\sigma$ is not
degenerate at points $\{s_0,\cdots,s_k \},$ i.e., the morphism
$\sigma^{n,0}: E^{n,0} \to E^{n,0} $ is injective at
$\{s_0,\cdots,s_k\}.$ Because $T$ is of dimension one and $\sigma$
is flat, we actually have
\begin{corollary}\label{flat-endomorphism-injective}
$\sigma$ is not degenerate at each point in $Y,$ i.e., the
morphism $\sigma^{n,0}: E^{n,0} \to E^{n-1,1}$ is injective
everywhere. Hence, the duality map $\sigma^{1,n-1}: E^{1,n-1}\to
E^{0,n} $ is surjective at each point of $Y.$ Furthermore,
$\sigma^{n+1}\equiv 0$ because $\sigma(E^{n,0})\equiv 0.$
\end{corollary}

\section{The Rigidity of Lefschetz Pencils of Calabi-Yau Varieties}

\begin{theorem}\label{zhangyi4}
Let $f:\ssX \rightarrow \P^1$ be a non-isotrivial Lefschetz pencil
of $n$-dimensional Calabi-Yau varieties. Then, the family $f$ must
be {\it rigid}.
\end{theorem}

\begin{proof}
Let $f^0:\ssX_0=f^{-1}(C_0) \rightarrow C_0=\P^1\setminus S$ be
the maximal smooth subfamily of $f.$ We have shown that $n$ must
be odd and there is a decomposition of the $\Q$-VHS
$$R_{prim}^nf^0_{*}(\Q)=\V\oplus (R_{prim}^nf^0_{*}(\Q))^{\pi_1(\P^1\setminus S)}.$$
with $(R_{prim}^nf^0_{*}(\Q))^{n,0}$ and
$(R_{prim}^nf^0_{*}(\Q))^{0,n}$ are in $\V_\C.$ Assume the
statement is not true,  we obtain a nontrivial deformation
$$
\begin{CD}
\ssX_0   \>\subset >> \sX \\
\V f^0 VV  \V V g V \\
C_0 \times \{0\}  \> \subset >> C_0 \times T^0
\end{CD}.
$$
where $T^0$ is a smooth quasi-projective curve. We then have a
nonzero flat $(-1,1)$-type endomorphism $\sigma$ of
$R_{prim}^nf^0_{*}(\C).$ As $R_{prim}^nf^0_{*}(\C)$ determines a
Higgs bundle $(E,\theta)$ on $C_0,$  we have:
\begin{lemma} \label{sigma-splitting}
Any nonzero flat $(-1,1)$-type endomorphism $\sigma \in \End(E)$
induces a splitting of the Higgs bundle
$$(E, \theta)=\Ker(\sigma)\bigoplus (\Ker(\sigma))^{\bot}.$$
The statement is also true if $C_0$ is replaced by a higher
dimensional quasi-projective variety $Y.$
\end{lemma}
\begin{proof}
The statement is a special case of the {\it generalized DSUY
correspondence} \ref{DSUY}. As the polarization $Q$ and the
endomorphism $\sigma$ both are flat under the Gauss-Manin
connection, we obtain a $\C$-splitting of the local system
$$R_{prim}^nf^0_{*}(\C)=\Ker(\sigma)\oplus (\Ker(\sigma))^{\bot}$$
which is compatible with the polarization $Q. $
$(\Ker(\sigma))^{\bot}$ is the orthogonal component of
$\Ker(\sigma)$ in $R_{prim}^nf^0_{*}(\C) ,$ i.e., the complete
reducibility. Therefore, if we regard $\sigma$ as a $\sO_C$-linear
map we have
$$R_{prim}^nf^0_{*}(\C)\otimes \sO_{C}=\Ker(\sigma)\oplus
(\Ker(\sigma))^{\bot}.$$ Restricting the Hodge filtration of the
VHS $R_{prim}^nf^0_{*}(\C)$ to these sub local systems and taking
the grading of the Hodge filtration, we actually have a
decomposition of the Higgs bundles. Moreover, it is really a
splitting of a {\it complex variation of Hodge structure} (cf.
\cite{Si2}).
\end{proof}

{\it Continue to prove the theorem.} $E^{0,n} \subset
\Ker(\sigma)$ along $M.$ Because of non-triviality of the
deformation of the family, there exists a point $s_0\in
\P^1\setminus S$ such that $\sigma: E^{n,0}|_{C_0} \rightarrow
E^{n-1,1}|_{C_0}$ is injective at $s_0$ by
\ref{flat-endomorphism-injective}. Thus,
$$E^{n,0}\nsubseteq
\Ker(\sigma) \ \ \mathrm{and} \ \ E^{0,n}\subset \Ker(\sigma)
$$ at $s_0.$ On the other hand,
$$\rho : \pi_1(\P^1\setminus S) \rightarrow  \Aut(\V_{\C,s_0})$$
is irreducible and both $E^{n,0}$ and $E^{0,n}$ are in $\V_\C$ by
\ref{zhangyi3}. We have a contradiction.
\end{proof}

Actually, we use two ideas in proving the theorem. The first idea
is that {\it a $(-1,1)$-type endomorphism of a complex polarized
VHS induces a splitting of the underlying $\C$-local system}, it
is implied in Zuo's theorem on the negativity of kernels of
Kodaira-Spencer maps of Hodge bundles (cf. \cite{Zuo0}). The
second idea is the well-known {\it Kazhdan-Margulis theorem} (cf.
\cite{Del1}): {\it For any Lefschetz pencil of odd dimensional
varieties, the image $\pi_1(\P^1\setminus S, s_0)$ is a Zariski
open set in $\mathrm{Sp}(\V_{\C,s_0}, (,))$ where $\V$ is the
$\Z$-local system of the vanishing cycles space.}

Altogether, the theorem \ref{zhangyi4} explores that there are
deep relations between two important objects in the algebraic
geometry: {\it rigidity} and {\it algebraic monodromy}. We start
to study these relations in \cite{ZY}.
\begin{definition}
Let $\V$ be a $\Q$-local system on a quasi-projective manifold $Y$
with monodromy representation
$$\rho: \pi_1(Y,y_0)\rightarrow \GL(V),\, \, V:=\V_{y_0}.$$
where $y_0$ is a fixed base point in $Y.$
\begin{myenumi}
\item The monodromy group $\pi_1(Y,y_0)^{mon}$ is defined to be
the Zariski closure of the smallest algebraic subgroup of $\GL(V)$
containing the monodromy representation $\rho(\pi_1(Y,y_0)). $
\item Assume that $V$ carries a nondegenerated bilinear form $Q$
which is symmetric (or anti-symmetric) and preserved by the
monodromy group.  We call the monodromy is {\it big} if the
connected component (including the identity) of
$\pi_1(Y,y_0)^{mon}$  acts irreducibly on $\V_\C. $
\end{myenumi}
\end{definition}

\begin{theorem}Let $f: X \to Y$ be a non-isotrivial smooth family of
Calabi-Yau $n$-folds.  Assume that $R_{prim}^nf_*(\Q)$ has a sub
$\Q$-VHS $\V$ with big monodromy and the bottom Hodge filtration
of the VHS $R_{prim}^nf_*(\Q)$ is in $\V.$ Then, the family $f$
must be rigid.
\end{theorem}
\section{The Rigidity of Strongly Degenerated Families}

\begin{definition}[Strongly degenerated families]
Let $f:\ssX \rightarrow C $ be a family over a smooth projective
curve $C$ with singular values $\{ c_0,\cdots,c_k\}.$ $f$ is
called {\it strongly degenerate} at $c_0$ if $f$ satisfies
following conditions:
\begin{myenumi}
    \item $X_i=f^{-1}(c_i)=X_{i1}+\cdots+X_{ir_i}$ is reduced and it is a union
    of transversally crossing smooth divisors for all $i,$
    i.e., $f$ is semistable.
    \item The cohomology of each component of the singular fiber $X_0=f^{-1}(c_0)$
    has pure type $(p,p).$
\end{myenumi}
\end{definition}

Let $f: \ssX \to C$ be a {\it strongly degenerated} family of
$n$-folds. Set $C_0= C \setminus \{ c_0,\cdots,c_k \},$ we have
the maximal smooth subfamily $f^0: \ssX_0 \to C_0$ of $f.$ Define
the fiber product family by $\pi: \sY\triangleq\ssX\times_C
\ssX\rightarrow C.$ $\pi$ is only degenerate at $\{c_0,\cdots,c_k
\},$ and $\pi^0: \sY_0=\ssX\times_{C_0}\ssX \to \ C_0$ is the
maximal smooth subfamily of $\pi.$ By the K\"unneth formula, we
obtain that
\begin{corollary}\label{product}
Let  $f:\ssX \rightarrow C$ be a family {\it strongly degenerate}
at $c_0 \in C.$ Then, the fiber product family $\pi$ is also {\it
strongly degenerate} at $c_0.$
\end{corollary}

Now, we study endomorphisms of Higgs bundles over a
quasi-projective variety $M.$ Let $(E,\theta)$ be a Higgs bundle
with positive Hermitian metric $H.$ Under the condition that
$(E,\theta)$ carries $\R$-structure, we describe $\End(E)^{-1, 1}$
precisely in the following way. As in section
\ref{product-in-the-section}, $\End(E)$ is also a Higgs bundle:
$$(\End(E)=\bigoplus_{r+s=0}\End(E)^{r,s},\theta^{end})=
(\bigoplus_{(p-p')+(q-q')=0}E^{p,q}\otimes E^{\vee
-p',-q'},\theta^{end})$$ with $ \theta^{end}(u\otimes
v^{\vee})=\theta(v)\otimes v^{\vee}+u\otimes
\theta_{\vee}(v^{\vee}).$ Thus,
\begin{equation}\label{(-1,1)}
\End(E)^{-1,1}=\bigoplus_{p+q=n}E^{p,q}\otimes E^{q-1,p+1},
\end{equation}
and $\End(E)^{-1,1} \subset  (E\otimes E)^{n-1,n+1}$ as vector
spaces.
\begin{lemma} \label{Kunneth-zhang}
Let $(E,\theta)$ be the Higgs bundle associated to
$R_{prim}^nf^0_{*}(\C).$ Then, by the K\"unneth formula
$H^*(X,\C)\otimes H^*(Z,\C)=H^*(X\times Z,\C),$ we have the
inclusion
\begin{equation}\label{Kunneth-formular} R_{prim}^nf^0_{*}(\C)
\otimes R_{prim}^nf^0_{*}(\C) \subset R^{2n}\pi^0_{*}(\C),
\end{equation}
which is compatible with Hodge structures. Therefore, as vector
spaces
\begin{equation}\label{Kunneth1}
\End(E)^{-1,1}\subset R^{n+1}\pi^0_{*}(\Omega_{\sY_0/C_0}^{n-1}).
\end{equation}
\end{lemma}

\begin{theorem}\label{main2-zhangyi}
Any non-isotrivial {\it strongly degenerated} family (not only for
families of Calabi-Yau varieties) must be rigid.
\end{theorem}
\begin{proof}
Let $f:\ssX\rightarrow C$ be a {\it strongly degenerated} family
of $n$-folds smooth over $C_0=C-\{c_0,\cdots,c_k \}.$ We have the
maximal smooth subfamily $f^0$ as before and $\pi,\pi^0$
respectively. Suppose that $f^0:\ssX_0\rightarrow C_0$ is
nonrigid, we then have a nonzero $(-1,1)$-type global section
$\sigma$ of $\End(R_{prim}^nf^0_{*}(\C)).$ By \ref{Kunneth-zhang},
$$0\neq \sigma\in
(R^{n+1}\pi^0_{*}(\Omega_{\sY_0/C_0}^{n-1}))^{\pi_1(C_0)}. $$ On the
other hand, there is a commutative diagram due to \ref{deligne}
$$
\begin{TriCDV}
{H^n(\sY,\C)}{\> i^*
>>}{H^n(\sY_0,\C)} {\SE  \overline{i}_{t}^* EE}{\SW  W i_{t}^* W}
{H^n(\sY_{t},\C)^{\pi_1(C_0,t)}}
\end{TriCDV}
$$
where $i_{t}:\sY_t\hookrightarrow \sY_0,$
$\overline{i}_{t}:\sY_t\hookrightarrow \sY$ are natural embedding;
and  $$H^n(\sY_{t},\C)^{\pi_1(C_0,t)}\cong
H^0(C_0,R^n\pi^0_{*}(\C))$$ is an isomorphism of Hodge structure.
For each $t\in C_0,$ as $\overline{i}^*_t$ is a surjective Hodge
morphism, we have the restriction maps
\begin{equation}\label{(p,q)-restriction}
r^{p,q}_t:H^q(\sY,\Omega^{p}_{\sY})\hookrightarrow
H^n(\sY,\C) \> \overline{i}_t^* >>
H^n(\sY_t,\C)\rightarrow H^q(\sY_t,\Omega^{p}_{\sY_t})
\end{equation}
where $p+q=2n.$ Thus, we obtain:
\begin{center}{\it The $(p,q)$-component group
$H^n(\sY_{t},\C)^{\pi_1(C_0,t)}$ is the image of
$H^q(\sY,\Omega^p_{\sY})$ under $r_t^{p,q}.
$}
\end{center}

Let $B\subset \sY$ be a $2n$-dimensional reduced nonsingular
algebraic cycle. Then,
$$\int_{B}\alpha\wedge\overline{\alpha}=\int_B\alpha|_{B}\wedge\overline{\alpha}|_{B} \, \forall \alpha\in
H^q(\sY,\Omega^{p}_{\sY})$$ and that
$\int_{B}\alpha\wedge\overline{\alpha}=0$ is equivalent to
$\alpha|_B=0. $ For all $\alpha\in H^q(\sY,\Omega^{p}_{\sY}),$ we
have $$ \alpha |_{D_j}=0 \, \forall j \Longleftrightarrow
\int_{Y_0}\alpha\wedge\overline{\alpha}=0$$ where
$Y_0=\pi^{-1}(c_0)=\sum D_j$ is a semistable singular fiber. On
the other hand, because all smooth fibers are homological
equivalent to $Y_0$ we have
$$\int_{Y_0}(\alpha\wedge\overline{\alpha})|_{Y_0}=\int_{\sY_t}(\alpha\wedge\overline{\alpha})|_{\sY_t}
\ \forall t\in C_0.$$

Hence, if $H^{p,q}(D_j)=0 \ \forall j$ then $\alpha |_{\sY_t}=0 \
\forall \alpha\in H^q(\sY,\Omega^{p}_{\sY}),$ i.e., the
restriction map
$$r^{p,q}_t:H^q(\sY,\Omega^{p}_{\sY})\rightarrow
H^q(\sY_t,\Omega^{p}_{\sY_t})\ \forall t\in C_0$$ is a zero map.

Now, $\pi$ is {\it strongly degenerate} at $c_0,$  the cohomology
of each component of $Y_0$ has only pure Hodge type. In
particular, $$H^{n-1,n+1}(D_j)=0\, \forall j.$$ Therefore, the
restriction map
$$r^{n-1,n+1}_t:H^{n+1}(\sY,\Omega^{n-1}_{\sY})\rightarrow
H^{n+1}(\sY_t,\Omega^{n-1}_{\sY_t}) \, \forall t \in C_0$$ has to
be zero, i.e.,
$$(R^{n+1}\pi^0_{*}(\Omega_{\sY_0/C_0}^{n-1}))^{\pi_1(C_0)}=0.$$
\end{proof}

From the proof the theorem, we actually obtain a stronger
property: {\it A non-trivial deformation of a family gives arise
to a non-trivial Betti-cohomology class in the cohomology group of
the total space of the self-product of the original family.} This
class is very special. It is not of Hodge type and does not vanish
along fibers. Thus, we have shown that there are relations between
the geometry of total spaces of families and the {\it rigidity}.

\begin{corollary}[Weakly Arakelov theorem]
Let $f:\sX \rightarrow C $ be a semistable family over a smooth
projective curve $C$ with at least one singular fiber, which is a
normal crossing divisor and each of its components is dominated by
a projective space. Then, $f$ is rigid.
\end{corollary}

{\bf Remark.} Our result is weaker than the original Arakelov
theorem, but it is for arbitrary higher dimensional varieties.

\begin{example}Let $F$ be a smooth homogenous polynomial with degree
$d.$ Any one parameter family in $\P^n$ of type
$$F(X_0,\cdots,X_n)+ t\prod_{i=1}^d X_{\tau (i)}=0 $$
is rigid, where $\tau: \{1,\cdots, d\} \to \{1,\cdots, n\}$ is an
injective map.
\end{example}
\section{Yukawa Couplings and Rigidity}
\subsection*{Yukawa couplings}
It is better to understand more about the endomorphism $\sigma$ in
section \ref{infi-rigid} which has deep background in string
theory. Let $\sigma^l$ be the $l$-iterated of $\sigma$-operator on
$E,$ then $\sigma^{l}\equiv 0$ for $l >> 0.$
\begin{proposition}\label{nonrigid}
Let $f^0:\sX_0\to Y_0$ be a smooth family of polarized Calabi-Yau
$n$-folds satisfying the $(\ast \ast)$-condition. Then, if $f$ is
nonrigid $$\sigma^n\equiv 0.$$
\end{proposition}
\begin{proof}
Otherwise, we have a nonzero flat endomorphism $\sigma^n.$ Denote
$ L:=f^0_*\Omega^n_{\sX_0/Y_0}.$ The endomorphism $\sigma^n$ in
fact is a global holomorphic section of the line bundle
$(L^*)^{\otimes 2}$ (under the holomorphic structure
$\overline{\partial}_E$ of the induced Higgs bundle), then
$(L^*)^{\otimes 2}$ is trivial. But it is impossible by the
results of \ref{zhangyi1} and \ref{WP-VHS}.
\end{proof}

\begin{definition}
Let $g :\sX\rightarrow Y$ be a smooth family of Calabi-Yau
$n$-folds over a quasi-projective manifold $Y$ and $(E,\theta)$ be
the Higgs bundle induced from the VHS $R_{prim}^ng_*(\C).$ The
{\it Yukawa coupling} is just the $n$-iterated Higgs field
$\theta^{n} : E\to E\otimes \Sym^n\Omega^1_Y.$
\end{definition}

The definition is a little different from classic literatures, but
they are compatible with each other. As the Higgs bundle
$(E,\theta)$ can be splitting into
$$(E,\theta)=(\bigoplus_{p+q=n}E^{p,q},\bigoplus \theta^{p,q})$$
with $\theta^{p,q}: E^{p,q}\to E^{p-1,q+1}\otimes \Omega_Y^1,$ the
$n$-iterated Higgs field $$\theta^{n,0}\circ\cdots\circ
\theta^{0,n}:E^{n,0}\>>>  E^{0,n}\otimes(\otimes^n \Omega_Y^1)$$
can factor through $\theta^n: E^{n,0}\>>>E^{0,n}\otimes
\Sym^n\Omega_Y^1 $ because $\theta\wedge \theta=0.$ So we can
formulate the {\it Yukawa coupling} as
$$\theta^n: \Sym^n\Theta_Y \rightarrow \sH om(E^{n,0},E^{0,n})=((R^0g_*\Omega^n_{\sX/Y})^*)^{\otimes 2}.$$
Assume that  $Y$ has a smooth compactification $\overline{Y}$ such
that $D_{\infty}=\overline{Y}\setminus Y$ is a normal crossing
divisor and the Higgs bundle $(E,\theta)$ over $Y$ has regular
singularities at $D_{\infty}.$ Then, $(E,\theta)$ is algebraic and
the {\it Yukawa coupling} is a global algebraic section of
$$((R^0g_*(\Omega^1_{\sX/Y}))^*)^{\otimes 2}\otimes
\Sym^n\Omega^1_{Y}. $$

Let $f^0:\ssX_0 \rightarrow C_0$ be a non-isotrivial smooth family
over a smooth quasi-projective curve  $C_0=C\setminus
\{s_1,\cdots,s_k \}.$  Suppose that $f^0$ is nonrigid. Then, we
have a deformation $F: \sX \rightarrow C_0\times T $ of $f^0$
where $\sX|_{\{C_0\}\times 0 }= \ssX_0,$ and a nonzero $\sigma$ on
the induced Higgs bundle over $C^0.$ By \ref{CY-sigma}, all
restricted families $$F_s: \sX|_{\{s\}\times T }\rightarrow
 \{s\}\times T \ \forall s\in C_0
$$ satisfy the $(\ast\ast)$-condition, and by \ref{nonrigid} there
are endomorphisms $$\sigma'_s : E|_{\{s\}\times T} \rightarrow
E|_{\{s\}\times T} \mbox{ (maybe zero)}$$ with
$(\sigma'_s)^n\equiv 0 \ \forall s\in C_0.$  On the other hand, we
have Higgs maps
$$\theta_t : E|_{C_0\times \{t\}} \rightarrow E|_{C_0\times \{t\}} \otimes \Omega^1_{C_0}$$
along $C_0.$ Because $\sigma'_{s}|_t$ comes from $\theta_{t}|_s$
for each $(s,t),$ we also can obtain $\sigma'_s $ by fixing $s\in
C_0 $ and varying $t \in T.$  Altogether, we have another proof of
the criterion of Viehweg-Zuo and Liu-Todorov-Yau-Zuo from the
construction of $\sigma.$
\begin{proposition}[\cite{VZ2,VZ4},\cite{LTYZ}]\label{yukawa-rigidity}
Let $f^0:\ssX_0\rightarrow C_0$ be a smooth family of Calabi-Yau
$n$-folds. If the {\it Yukawa coupling} is not zero at some points,
then the family $f^0$ must be rigid.
\end{proposition}

{\bf Remarks.} That the {\it Yukawa coupling} is nonzero at some
points implies that the family $f$ satisfies the $(\ast
\ast)$-condition, so $f$ is automatically non-isotrivial. Denote
$\sZ=\{s\in C_0 \ | \ \theta^n_s=0\}.$ $\sZ$ is either a set of
finite points or the total $C_0.$  The criterion actually says
that if $\sZ$ is a finite set then $f$ is rigid. We shall point
out that the converse statement is not always true: Using the
covering trick, Viehweg-Zuo recently constructed some rigid
families of Calabi-Yau varieties with the {\it Yukawa coupling}
identifying with zero (cf. \cite{VZ7}).

We would like to introduce a stronger result of Viehweg-Zuo, they
dealt with more general cases: not only for Calabi-Yau manifolds
but also for projective manifolds with semi-ample canonical line
bundle or minimal models of general type.
\begin{criterion}\label{VZ-criterion}\cite[Corollary 6.5]{VZ2} or \cite[Corollary
8.4]{VZ4}. Let $h$ be a fixed Hilbert polynomial of degree $n$ and
$\ssM_h$ be a coarse moduli space of polarized $n$-folds with
semi-ample canonical line bundle (including Calabi-Yau manifolds)
or of general type. Assume $\ssM_h$ has a nice compactification
and carries a universal family
$\pi:\ssX\to\ssM_h $ 
 (In the real situation, one needs to work on stacks). Let
$k_{\ssM}$ be the largest integer such that the $k_{\ssM_h}$-times
iterated Kodaira-Spencer map of this universal family is not zero
(obviously, $1\leq k_{\ssM_h}\leq n$). If the $k_{\ssM_h}$-times
iterated Kodaira-Spencer map for a family $f:X\to Y$ is not zero,
then the family $f$ must be rigid.
\end{criterion}

For higher dimensional base,  one always consider the following
condition (good partial compactification cf. \cite{VZ2}):  Let $U$
be a manifold and $Y$ be a smooth projective compactification of
$U$ with a reduced normal crossing divisor boundary
$D_{\infty}=Y\setminus U.$ Starting with a smooth family $f: V\to
U,$ we first choose a smooth projective compactification $X$ of
$V,$ such that $f: V\to U$ extends to $f:X\to Y.$ Then, one leaves
out codimension $2$ subschemes of $Y$ such that the restricted
morphism is flat (cf. \cite{VZ2,VZ4}).

A theorem of Landman-Katz-Borel says that local monodromies for
the VHS $R^nf_*{\Q_V}$ around $D_\infty$ are all quasi-unipotent
(cf. \cite{Sch}), so we have Deligne's quasi-canonical extension
of the VHS, i.e., that the real part of the eigenvalues of the
residues around the components of $S$ lies in $[0,1)).$ Take
grading of the filtration, we get the quasi-canonical extension of
the Higgs bundle
$$(\bigoplus_{p+q=n} \overline{E}^{p,q}, \  \bigoplus \overline{\theta}^{p,q})$$
where $\overline{\theta}^{p,q}: \overline{E}^{p,q} \to
\overline{E}^{p-1, q+1}\otimes \Omega^1_{Y}(\log S). $ We still
call  $$\overline{\theta}^n : \overline{E}^{n,0} \>>>
\overline{E}^{0,n}\otimes \Sym^n\Omega^1_Y(\log S)$$ the {\it
Yukawa coupling}.  If there is no confusion, we also formulate the
{\it Yukawa coupling} as $$\overline{\theta}^n : \Sym^n (T_Y(-\log
S)) \to \overline{E}^{0,n}\otimes (\overline{E}^{n,0})^\vee. $$
One has $\overline{E}^{p,q}=R^qf_*\Omega^p_{\sX/Y}(\log \Delta)$
and
$$\Omega^n_{\sX/Y}(\log \Delta) = \omega_{\sX/Y}(\Delta_{\rm red}
- \Delta)$$ where $\Delta=f^* S$ (cf. \cite{Gr}). Denote
$\sL:=\Omega^n_{\sX/Y}(\log \Delta).$ The {\it Yukawa coupling} is
$$ \overline{\theta}^n: R^0f_*\Omega^n_{\sX/Y}(\log \Delta)\>>> R^nf_*{\sO_{\sX}} \otimes \Sym^n\Omega^1_Y(\log S), $$
or it is $$\overline{\theta}^n: \Sym^n T^1_Y(-\log S) \>>>
R^nf_*{\sO_{\sX}} \otimes (f_*\sL)^{-1}. $$ If $\Delta$ is reduce,
then $\sL=\omega_{X/Y}$ and the {\it Yukawa coupling} then is
$$ \overline{\theta}^n: \Sym^n T^1_Y(-\log S) \>>>
R^nf_*{\sO_{\sX}} \otimes (f_*\omega_{X/Y})^{-1}. $$

We have a generalization of \ref{yukawa-rigidity} from the
Viehweg-Zuo criterion.
\begin{theorem}\label{VZ-rigid}Under assumptions made above, let $f:\sX \to Y$ be a family of
$n$-dimensional projective varieties with a general fiber having
semi-ample canonical line bundle or being of general type. If the
{\it Yukawa coupling} $\overline{\theta}^n(f)$ is not zero, then
the family $f$ must be rigid.
\end{theorem}
\begin{proof}
First, we have the tautological sequence
\begin{equation*}\label{log-complex}
0 \>>> f^*\Omega^1_Y(\log S) \>>> \Omega^1_\sX(\log \Delta) \>>>
\Omega^1_{\sX/Y}(\log \Delta) \>>> 0,
\end{equation*}
and the wedge product sequences
\begin{multline}\label{exact}
0\>>> {f}^*\Omega^1_Y(\log S)\otimes \Omega^{p-1}_{\sX/Y}(\log
\Delta) \>>> \\ {\mathfrak g \mathfrak r}(\Omega_\sX^p(\log
\Delta)) \>>> \Omega_{\sX/Y}^p(\log \Delta)\>>> 0
\end{multline}
where $${\mathfrak g \mathfrak r}(\Omega_\sX^p(\log \Delta))=
\Omega_\sX^p(\log \Delta) /f^*\Omega^2_Y(\log S)\otimes
\Omega^{p-2}_{\sX/Y}(\log \Delta).$$  We study various sheaves
$$F^{p,q}=R^qf_*(\Omega^{p}_{\sX/Y}(\log \Delta)\otimes\sL^{-1})$$
for $\sL=\Omega^n_{\sX/Y}(\log \Delta),$ together with edge
morphisms $$\tau_{p,q}:F^{p,q}\to F^{p-1,q+1}\otimes
\Omega^1_{Y}(\log S) \,\forall (p,q) \mbox{ with } p+q=n$$ induced
by the exact sequences (\ref{exact}), tensored with $\sL^{-1}.$
Explained in \cite[4.4 iii]{VZ2},
$$R^{q}f_*(\wedge^{n-p}T_{\sX/Y}(-\log
\Delta))=R^qf_*(\Omega^{p}_{\sX/Y}(\log \Delta)\otimes \sL^{-1})$$
for $\forall (p,q) \mbox{ with } p+q=n,$ and all $F^{p,q}$ are
indeed from the deformations of the family $f:\sX \to Y.$

Over $U=Y\setminus S,$ edge morphisms $\tau_{p,q}$ also can be
obtained in the following way. Consider the exact sequence $$ 0
\to T_{V/U} \to T_V \to f^*T_U \to 0.$$ It induces sequences
$$
0 \>>> \bigwedge^{n-p+1}T_{V/U} \>>> \tilde T^{n-p+1}_V \>>>
\bigwedge^{n-p}T_{V/U}\otimes f^*T_U \>>> 0,
$$
where $\tilde{T}^{n-p+1}_{V}$ is a subsheaf of
$\bigwedge^{n-p+1}T_V$ for each $p.$ The edge morphisms are $$
\tau^\vee_{p,q}:(R^{q}f_*(\wedge^{n-p}T_{V/U})) \otimes T_U \to
R^{q+1}f_*(\wedge^{n-p+1}T_{V/U}),$$ they are the wedge product
with the Kodaira-Spencer class. Tensoring with $\Omega_U^1,$ we
get back $\tau_{p,q}|_U.$ Thus, $$\tau_{p,q} : F^{p,q} \to
F^{p-1,q+1}\otimes \Omega_Y^1(\log S) \,\forall (p,q)$$ are the
{\it log Kodaira-Spencer maps}, and
\begin{multline*}
\tau^m : F^{n,0}=\sO_Y \> \tau_{n,0} >> F^{n-1,1}\otimes
\Omega_{Y}^1(\log S)  \\ \> \tau_{n-1,1} >>  F^{n-2,2}\otimes
\Sym^2(\Omega^1_Y(\log S))\>>> \\
\cdots  \> \tau_{n-m+1,m-1}>> F^{n-m,m}\otimes
\Sym^m(\Omega^1_Y(\log S))
\end{multline*}
is the {\it $m$-times iterated log Kodaira-Spencer class.}
Finally, we have a factor map by
$$
\begin{TriCDA}
{\Sym^n(T_Y(-\log S))} {\SW\tau^n(f) WW}{\SE E
\overline{\theta}^n(f) E} {R^n f_*(\sL^{-1})}{\>>>}
{R^nf_*{\sO_{\sX}} \otimes (f_*\sL)^{-1}}
\end{TriCDA}
$$
and $$\overline{\theta}^{n}(f)\neq 0 \Longrightarrow \tau^n(f)\neq
0 \Longrightarrow {k_{\ssM}=n} \Longrightarrow
\tau^{k_{\ssM}}_{\ssM}\neq 0.$$ Hence, the family $f$ must be
rigid.
\end{proof}

{\bf Remark.} Except if $\rank \overline{E}^{n,0}=1,$ it is
difficult to prove a similar result only by regarding variations
of Hodge structures. So one should study geometric deformations
more precisely.
\begin{example}\label{ltyz}
Let $F,$ $H$ be two homogenous polynomials of degree $n+2$ such
that $F$ defines a nonsingular hypersurface in $\P^{n+1}$ and $H$
is not in the Jacobian ideal of $F.$ Consider a special family
$F(t)=F+tH, t\in \P^1 $ which satisfies that $H^n$ is not in the
Jacobian ideal of $F+\mu H$ for some $\mu \in \P^1.$ Then, the{\it
Yukawa coupling} at $\mu\in\P^1$ is not zero (cf. \cite{LTYZ}),
and the family is rigid.
\end{example}

\subsection*{Residues of Higgs fields}

Let $f:\ssX\rightarrow C $ be a family smooth over $C_0.$ Let
$(\sV,\nabla)$ be a locally free sheaf with the Gauss-Manin
connection given by the weight $n$ polarized VHS
$R_{prim}^nf_*(\Q_{X_0})$ where $\ssX_0=f^{-1}(C_0).$ As all
monodromies for the VHS are quasi-unipotent, one has Deligne's
quasi-canonical extension $(\overline{\sV},\overline{\nabla})$
over $C$ with $\overline{\nabla} : \overline{\sV} \rightarrow
\overline{\sV}\otimes \Omega_{C}^1(\log S).$ $(\sV,\nabla)$
corresponds to a regular filtered Higgs bundle $\{
(E,\theta)_{\alpha}\}.$ For convenience, we restrict $f$ to a unit
disk $\Delta$ and study the local degenerated family
$f_{loc}:\ssX_{loc} \to (\Delta,t)$ which is smooth over
$\Delta^*=\Delta-\{0\}.$ $\overline{\sV}|_0$ can be represented by
$\sV|_{\Delta^*}$ and the local monodromy $T$ extends naturally to
be an endomorphism on $\overline{E}.$ $T_0=T|_{\sV|_0}: \sV|_0 \to
\sV|_0$ is the restriction endomorphism (cf. \cite{Si2}).
\begin{definition}[Residues of $(\overline{\sV}, \overline{\nabla})$ \cite{Ste}]
For an integrable logarithmic connection $ \overline{\nabla} :
\overline{\sV} \to \overline{\sV} \otimes \Omega^{1}_{\Delta}
(\log0),$ one has a composite map
$$  (id_{\overline{\sV}} \otimes R_0)\circ \overline{\nabla} : \overline{\sV} \>
\overline{\nabla}
>> \overline{\sV} \otimes \Omega^{1}_{\Delta} (\log 0)  \>  id_{\overline{\sV}} \otimes R_0
>> \overline{\sV}|_0 \ ,$$
where $R_0:\Omega^1_{\Delta}(\log 0)\to\C$ is defined by
$R_0(hdt/t)=h(0).$ The composite map is zero at $t\overline{\sV},$
then it defines the {\it residue morphism}
$$\Res_0(\overline{\nabla}): \overline{\sV}|_0 \to
\overline{\sV}|_0.$$
\end{definition}
The local monodromy $T$ around $0$ is called {\it unipotent} if $
(T-1)^{k+1}=0$ and $(T-1)^{k}\neq 0$ for a fixed integer $k\in
[0,n].$ If $k=n, $ the monodromy $T$ is called  {\it maximal
unipotent} (cf. \cite{Del7}).

Suppose that $T$ is {\it unipotent.} Let $$N:=\log
T=\sum_{j=1}^{k}(1/j)(-1)^j(T-\Id)^j.$$ Then, one has $N^{k+1}=0$
and $N^{k}\neq 0.$ The canonical extension $(\overline{\sV},
\overline{\nabla})$ is generated  by all sections
$$\widetilde{v}(t)=\exp(\frac{-\log t}{2\pi\sqrt{-1}}N)\cdot v$$
where $v$ is a flat section (multiplicative values) of $\sV$ over
$\Delta^*.$ Moreover,
$$\widetilde{v}(te^{2\pi\sqrt{-1}})=\widetilde{v}(t).$$
\begin{lemma}\cite[Theorem II 3.11]{Del0}\label{residue}
Let $(\sV,\nabla)$ be a locally free sheaf with the Gauss-Manin
connection over $\Delta^*$ given by a polarized VHS with {\it
unipotent} monodromy $T.$ Then,
$$\Res_{0}(\overline{\nabla})=\frac{-1}{2\pi\sqrt{-1}} \log T= \frac{-1}{2\pi\sqrt{-1}} N.$$
\end{lemma}
\begin{definition}[Residues of Higgs fields \cite{Si2}]
Let $(V,H,D)$ be a {\it tame harmonic bundle} over $\Delta^*$ and
$\{(E,\theta)_{\alpha}\}$ be the induced regular filtered Higgs
sheaf on $\Delta.$ Denote $\overline{E}=\cup E_{\alpha}$ and
$\overline{\theta}=\cup \theta_{\alpha}.$ The composite map
$$  (id_{\overline{E}} \otimes R_0)\circ \overline{\theta} : \overline{E}
\> \overline{\theta} >> \overline{E}\otimes\Omega^{1}_{\Delta}
(\log 0)
  \> id_{\overline{E}} \otimes R_0 >> \overline{E}|_0 $$
induces the residue $\Res_0(\overline{\theta}): \overline{E}|_0
\to \overline{E}|_0. $ Precisely, the residue is from
$$0\>>> \Omega^1_{\Delta}\otimes \sE nd(\overline{E})\>>> \Omega^1_{\Delta}(\log 0)
\otimes  \sE nd(\overline{E}) \> R_0\otimes id >>  \sE
nd(\overline{E})|_{0} \>>>0.$$
\end{definition}

{\bf Remark.} If a Higgs bundle $(E,\theta)$ is induced from a
polarized VHS, $\Res_0(\overline{\theta})$ is automatically
nilpotent.
\begin{proposition}[Schmid-Simpson \cite{Sch}\cite{Si2}]\label{S-S}
Under the isomorphism $\overline{E}|_0 \cong \overline{\sV}|_0,$
the nilpotent part of $\Res_0(\overline{\nabla})$ is isomorphic to
the nilpotent part of $\Res_0(\overline{\theta}).$
\end{proposition}
The proof is simple if the local monodromy $T$ around $0$ of the
VHS $(\sV,\nabla)$ over $\Delta^*$ is {\it unipotent}: The
filtration of extensive Higgs bundles is simple and jumps only at
$\alpha=0,$ so one has $$\Res_0(\overline{\theta}) \cong
\Res_0(\overline{\nabla}) =\frac{-1}{2\pi\sqrt{-1}} \log T.$$
\begin{theorem}\label{max-uni}
Let $f:\ssX\rightarrow C$ be a family of $n$-dimensional
Calabi-Yau varieties. If $f$ admits a degeneration with {\it
maximal unipotent monodromy}, then the family $f$ must be rigid.
\end{theorem}
\begin{proof}
By \ref{residue} and \ref{S-S}, the {\it Yukawa coupling} is
nonzero if the local family $f_{loc}:\ssX_{loc} \to \Delta$ is
degenerate at $0$ with a {\it maximal unipotent monodromy} (cf.
\cite{ZY}).
\end{proof}

\begin{example}[cf. \cite{LTY}]\label{lty}
Lian-Todorov-Yau recently studied  a family of Calabi-Yau
varieties in $\mathbb{P}^{n+k} $ with $n \geq 4 , k \geq 1$
determined by
$$G_{1,t}=tF_{1}-\prod_{i=0}^{n_{1}}x_{i}=0,... ,
 G_{k,t}=tF_{k}-\prod _{j=n_{1}+..+n_{k-1}}^{n_{k}}x_{j}=0$$ for $t \in
\P^1$ and  $[x_0,\cdots, x_{n+k}]\in \P^{n+k},$ where the system
$$F_{1}=...=F_{k}=0 \mbox{ with $n_{i}=\deg F_{i}\geq2$ and $\sum n_{i}=n+k+1$}$$ defines a smooth Calabi-Yau variety.
This family must admit a degeneration with {\it maximal unipotent
monodromy}, hence it is a rigid family.
\end{example}
\begin{example}[cf. \cite{Del7}]\label{Morrsion} As a special case of Lian-Todorov-Yau, Morrison
showed that the {\it Yukawa coupling} is nonzero at $0$ for a
family defined by
$$X_0^5+X_1^5+X_2^5+X_3^5+X_4^5-5\lambda X_0 X_1 X_2 X_3 X_4=0, [X_0,X_1,X_2,X_3,X_4]\in \P^4$$
where $\lambda\in \P^1$ is a parameter. By \ref{lty}, this family
admits a degeneration with {\it maximal unipotent monodromy} at
$0.$ On the other hand, this family is rigid also by \ref{ltyz},
because $(X_0 X_1 X_2 X_3 X_4)^3$ is not in the Jacobian ideal of
$X_0^5+X_1^5+X_2^5+X_3^5+X_4^5.$
\end{example}

{\bf Remark.} The two families in examples \ref{lty} and
\ref{Morrsion} are all {\it strongly degenerated} (after a
semistable reduction), so that their rigidity also follows from \ref{main2-zhangyi}.\\

\section{Acknowledgements}



I am very grateful to Shing-Tung Yau and Kang Zuo for instruction,
encouragement and sharing ideas with me. I wish to thank Eckart
Viehweg for valuable suggestions on moduli spaces and Shafarevich
problems, Ngaiming Mok for pointing out that the non-positivity of
curvatures in the horizontal direction of period domains would
underly the validity of the criterion for rigidity. I also wish to
thank Andrey Todorov and Kefeng Liu for helpful discussions.

\bibliographystyle{plain}

\begin{thebibliography}{XXX}
\bibitem{Ara} S.Ju. Arakelov,  {\it Families of algebraic curves with fixed degeneracies} (Russian),
Izv. Akad. Nauk SSSR Ser. Mat. {\bf 35} (1971) 1269--1293,
MR0321933, Zbl 0248.14004.
\bibitem{BV} E. Bedulev \&  E. Viehweg, {\it On the Shafarevich conjecture
for surfaces of general type over function fields}, Inv. Math.
{\bf 139} (2000) 603-615, MR1738062.
\bibitem{CG} M. Cornalba \& P. Griffiths, {\it Analytic cycles and vector bundles on non-compact algebraic varieties},
Inv. Math. {\bf 28} (1975) 1--106,  MR0367263, Zbl 0293.32026.
\bibitem{Do1} S.K. Donaldson, {\it Anti self-dual Yang-Mills connections over complex algebraic surfaces
and stable vector bundles}, Proc. London Math. Soc. (3) {\bf
50(1)} (1985) 1--26, MR0765366, Zbl 0529.53018.
\bibitem{Del0} P. Deligne,  {\it \'Equations diff\'erentielles \'a points singuliers
r\'eguliers} (French), Lecture Notes in Mathematics, {\bf 163}
Springer-Verlag, Berlin-New York, 1970, MR0417174, Zbl 0244.14004.
\bibitem{Del}  P. Deligne,  {\it Th\'eorie de Hodge} II, I.H.\'E.S. Publ.Math.
{\bf 40} (1971) 5--57, MR0498551, Zbl 0219.14006.
\bibitem{Del1} P. Deligne, {\it  La conjecture de Weil}  I, I.H.\'E.S. Publ.Math. {\bf
43 } (1974) 273-307,   MR0387282, Zbl 0314.14007.
\bibitem{Del7} P. Deligne, {\it Local behavior of Hodge structures at
infinity}, in \textquoteleft Mirror symmetry II', 683--699, AMS/IP
Stud. Adv. Math., {\bf 1}, 1997,  MR1416353, Zbl 0939.14005.
\bibitem{EM} P. Eyssidieux \& N.M.  Mok, {\it Characterization of certain holomorphic geodesic cycles on
Hermitian locally symmetric manifolds of the noncompact type}, in
\textquoteleft Modern methods in complex analysis', 79-117, Ann.
of Math. Stud., {\bf 137}, Princeton Univ. Press, Princeton, NJ,
1995, MR1369135, Zbl 1013.32013.
\bibitem{E1}  P. Eyssidieux, {\it La caract\'eristique d'Euler du complexe de Gauss-Manin} (French),
J. Reine Angew. Math. {\bf 490} (1997) 155--212,  MR1468929, Zbl
0886.32013.
\bibitem{E2} P. Eyssidieux, {\it K\"ahler hyperbolicity and variations of Hodge
structures}, in \textquoteleft New trends in algebraic geometry
(Warwick, 1996)', 71--92, London Math. Soc. Lecture Note Ser.,
{\bf 264}, ed. K. Hulek, F. Catanese,  C. Peter, \& M. Reid
Cambridge Univ. Press, Cambridge, 1999, MR1714821, Zbl 0952.32014.

\bibitem{Falt}  G. Faltings,  {\it Arakelov's theorem for Abelian
varieties}, Inv. Math. {\bf 73} (1983) 337-347, MR0718934, Zbl
0588.14025.
\bibitem{Gr} P. Griffths,  {\it Topices in Transcendental Algebraic
Geometry}, Ann. of Math. Stud., {\bf 106} 1984 Princeton Univ.
Press. Princeton, N.J., Zbl 0528.00004.
\bibitem{Gr2} P. Griffiths,  {\it Periods of integrals on
algebraic manifolds} III,  I.H.\'E.S. Publ. Math. {\bf 38} (1970)
125-180,  MR0282990, Zbl 0212.53503.
\bibitem{Gro} A. Grothendieck  et al (with P. Deligne \& N. Katz)
{\it Seminaire de Geometrie Algebrique du Bois-Marie}, SGA 7 Parts
I and II, Springer Lecture Notes in Math. {\bf 288-340}, 1971 to
1977, MR0354656, Zbl 0237.00013 ; Zbl 0258.00005.
\bibitem{GM} M. Goresky \& R. MacPherson, {\it Stratified Morse Theory}, Springer,
1988, MR0932724,  Zbl 0639.14012.

\bibitem{JY} J. Jost  \&  S.-T. Yau, {\it Harmonic mappings and algebraic varieties over function fields},
Amer. J. Math. {\bf 115(6)} (1993) 1197--1227, MR1254732, Zbl
0824.14008.
\bibitem{JZ1}   J. Jost \&  K. Zuo,
{\it Harmonic maps of infinite energy and rigidity results for
representations of fundamental groups of quasiprojective
varieties}, J. Diff. Geom. {\bf 47(3)} (1997) 469--503, MR1617644,
Zbl 0911.58012.
\bibitem{JZ2} J. Jost \&  K. Zuo, {\it Harmonic maps into Bruhat-Tits buildings and
factorizations of $p$-adically unbounded representations of
$\pi\sb 1$ of algebraic varieties},  J. Alg. Geom. {\bf 9(1)}
(2000) 1--42, MR1713518, Zbl 0984.14011.
\bibitem{JZ3} J. Jost \& K. Zuo, {\it Arakelov type inequalities for
Hodge bundles over algebraic varieties, Part 1: Hodge bundles over
algebraic curves}, J. Alg. Geom. {\bf 11} (2002), 535-546,
MR1894937.
\bibitem{Ko} S.J. Kov\'acs, {\it  Families over base with
birationally nef tangent bundle}, J. Reine Angew. Math. {\bf 487}
(1997) 171-177,  MR1464907, Zbl 0922.14024
\bibitem{Kovacs} S.J. Kov\'acs, {\it Logarithmic vanishing theorems
and Arakelov-Parshin boundedness for singular varieties},
Compositio Math. {\bf 131(3)} (2002) 291--317, MR1905025, Zbl
pre01764836.
\bibitem{LTY} B.H. Lian,  A. Todorov \& S.-T. Yau, {\it Maximal unipotent monodromy for
 complete intersection CY manifolds},  preprint 2000,
 Math.AG/0008061.
\bibitem{LTYZ} K. Liu, A. Todorov,  S.-T. Yau \&  K. Zuo,
{\it Shafarevich's conjecture for CY Manifolds} I,  preprint 2003,
Math.AG/0308209.
\bibitem{Mi} Li. Migliorini, {\it A smooth family of minimal surfaces of general type over a curve of genus at Most one is trivial},
J. Alg. Geom. {\bf 4} (1995) 353-361, MR1311355, Zbl 0834.14021.
\bibitem{Mok3}  N.M. Mok, {\it Characterization of certain holomorphic geodesic cycles on quotients
of bounded symmetric domains in terms of tangent subspaces},
Compositio Math. {\bf 132(3)} (2002) 289--309,  MR1918134, Zbl
1013.32013.
\bibitem{NS} M.S. Narasimhan \& C.S. Seshadri, {\it Stable and unitary vector bundles on a compact
Riemann surface}, Ann. of Math. (2), {\bf 82} (1965) 540--567,
MR0184252, Zbl 0171.04803.

\bibitem{O-V} K. Oguiso \& E. Viehweg, {\it On the isotriviality of families of elliptic
surfaces} J. Alg. Geom. {\bf 10(3)} (2001) 569--598, MR1832333.
\bibitem{Par} A.N. Par\v sin, {\it Algebraic curves over function fields} I (Russian),
Izv. Akad. Nauk SSSR Ser. Mat. {\bf 32} (1968) 1191--1219,
MR0257086, Zbl 0181.23902.
\bibitem{Pe}  C. Peters, {\it Rigidity for variations of Hodge
structures and Arakelov-type finiteness theorems}, Compositio
Math. {\bf 75} (1990) 113-126,  MR1059957, Zbl 0743.14006.
\bibitem{Sch} W. Schmid, {\it Variation of Hodge structure : The singularities
of the period mapping}, Inv. Math. {\bf 22} (1973) 211--319,
MR0382272, Zbl 0278.14003.
\bibitem{Sh} I.R.  Shafarevich, {\it Collected Mathematical Papers},
New York, Springer-Verlag 1989, MR0977275, Zbl 0669.12001.
\bibitem{Si1}  C. Simpson, {\it Constructing variations of Hodge structure using Yang-Mills theory
and applications to uniformization}, Journal of the AMS {\bf 1}
(1988) 867-918,  MR0944577, Zbl 0669.58008.
\bibitem{Si2} C. Simpson, {\it Harmonic bundles on noncompact
curves}, Journal of the AMS {\bf 3} (1990) 713--770, MR1040197,
Zbl 0713.58012.
\bibitem{Si3}  C. Simpson, {\it Higgs bundles and local
system}, I.H.\'E.S. Publ.Math. {\bf 75} (1992) 5-95,  MR1179076,
Zbl 0814.32003.

\bibitem{Ste} J. Steenbrink,  {\it Limits of Hodge structures},
Inv. Math. {\bf 31(3)} (1975/76) 229--257,   MR0429885, Zbl
0303.14002.
\bibitem{SZ} H.M.  Saito \&  S.   Zucker,  {\it Classification of non-rigid
families of K3 surfaces and finiteness theorem of Arakelov type},
Math. Ann. {\bf 298} (1993) 1-31,  MR1087233, Zbl 0697.14024.
\bibitem{Ti} G. Tian, {\it Smoothness of the universal deformation space
of Calabi-Yau manifolds and its Petersson-Weil metric}, in
\textquoteleft Math. Aspects of String Theory', ed. S.-T. Yau,
World Scientific 1998, 629-346, MR0915841, Zbl 0696.53040.
\bibitem{To89} A.  Todorov,  {\it The Weil-Petersson geometry of the moduli
space of $SU(n)$ $\geq3$ (Calabi-Yau Manifolds)} I, Comm. Math.
Phys. {\bf 126} (1989) 325-346, MR1027500, Zbl 0688.53030.

\bibitem{UY} K. Uhlenbeck \& S.-T. Yau, {\it On the existence of Hermitian-Yang-Mills connections in stable vector
bundles}, in \textquoteleft Frontiers of the mathematical sciences
1985' (New York, 1985), Comm. Pure Appl. Math. {\bf 39(S)} (1986)
S257--S293, MR0861491, Zbl 0615.58045.
\bibitem{V95} E. Viehweg, {\it Quasi-Projective Moduli for
Polarized Manifolds}, Ergebnisse der Mathematik, 3. Folge {\bf 30}
(1995), Springer Verlag, Berlin-Heidelberg-New York, MR1368632,
Zbl 0844.14004.
\bibitem{VZ1} E. Viehweg \&  K. Zuo, {\it On the isotriviality of families
of projective manifolds over curves}, J. Alg. Geom. {\bf 10}
(2001), 781-799,  MR1838979.
\bibitem{VZ2} E. Viehweg \&  K. Zuo, {\it Base spaces of non-isotrivial families of smooth minimal
models}, in \textquoteleft Complex geometry' (G\"ottingen, 2000),
279--328, Springer, Berlin, 2000,  MR1922109, Zbl 1006.14004.
\bibitem{VZ3}E. Viehweg \&  K. Zuo, {\it On the Brody hyperbolicity of
moduli spaces for canonically polarized manifolds} Duke Math
Journal {\bf 118(1)} (2003) 103-150, MR1978884.
\bibitem{VZ4} E. Viehweg \&  K. Zuo, {\it Discreteness of minimal models of Kodaira dimension zero
and subvarieties of moduli stacks}, Surveys in differential
geometry, Vol. VIII (Boston, MA, 2002), 337--356, Surv. Diff.
Geom., {\bf VIII}, Int. Press, Somerville, MA,
 2003, MR2039995.
\bibitem{VZ6} E. Viehweg \&  K. Zuo, {\it Families over curves with a strictly maximal Higgs
field}, Asian Journal of Mathematics, {\bf 7(4)} (2003) 575-598,
MR2074892.
\bibitem{VZ7}E. Viehweg \&  K. Zuo, {\it Complex multiplication, Griffiths-Yukawa couplings, and
rigidity for families of hypersurfaces}, preprint 2003,
Math.AG/0307398.
\bibitem{ZY} Y. Zhang,  {\it On families of Calabi-Yau manifolds}, Ph.D. Thesis, 2003, The Chinese University of
Hong Kong.
\bibitem{Z} Q. Zhang,  {\it Global holomorphic one-forms on projective manifolds
with ample canonical bundles}, J. Alg. Geom. {\bf 6} (1997)
777-787, MR1487236, Zbl 0922.14008.
\bibitem{Zuo1} K. Zuo, {\it Representations of fundamental groups of
algebraic varieties}, Lecture Notes in Mathematics, {\bf 1708},
Springer-Verlag, Berlin, 1999, MR1738433, Zbl 0987.14014.
\bibitem{Zuo0}  K. Zuo, {\it  On the negativity of kernels of
Kodaira-Spencer maps on Hodge bundles and applications}, Asian
Journal of Mathematics. {\bf 4(1)} (2000) 279-302, MR1803724, Zbl
0983.32020.
\bibitem{Zuo5} K. Zuo, {\it On families of projective manifolds} ({\it joint works with Eckart
Viehweg}), Manisctipt, a full time talk at France-Hong Kong
Geometry Conference, Feb 18-Feb 22, 2002, Hong Kong.
\end{thebibliography}

\end{document}